\documentclass[12pt]{amsart}
\usepackage{CFmath-english}
\usepackage{longtable}

\newcommand{\Perfect}{\operatorname{\widetilde P}}
\newcommand{\NIPer}{\operatorname{P}}
\newcommand{\NIPerM}{\operatorname{N}}

\title{Counting periodic points of rational self-maps of projective varieties}
\author{Xiao LI}
\address{Centre de Math\'ematiques Laurent Schwartz,
  \'Ecole Polytechnique,
  91128 Palaiseau Cedex, France}
\email{\href{mailto:xiao.li.cmls@polytechnique.edu}{xiao.li.cmls@polytechnique.edu}}
\date{November 24, 2026}

\begin{document}

\begin{abstract}
  We study the asymptotic behaviour of the number of isolated periodic points of a dominant rational
  self-map of a smooth projective variety. We first give an upper bound for this number in terms of
  the degrees of the iterates of the map, thereby extending a theorem of Dinh--Nguyen--Truong to
  arbitrary characteristic. Our approach is purely algebraic and bypasses the theory of densities of
  positive closed currents.

  We then establish the exact asymptotic growth of isolated periodic points of birational self-maps
  of surfaces over a field of characteristic 0. To do so, we extend Jonsson--Reschke's argument to
  any finitely generated field over $\mathbb Q$ using Moriwaki's theory of heights.
\end{abstract}

\maketitle

\section{Introduction}

Counting periodic points is a basic problem in dynamical systems. It gives
a quantitative invariant of the complexity of the dynamics and is closely
related to entropy and to the growth of cohomological data.

In this paper, we study this question in an algebraic framework. Let $X$ be a normal projective
variety of dimension $m$ defined over an algebraically closed field $K$, and let
$f\colon X\dasharrow X$ be a dominant rational self-map. Recall that $f$ is determined by its graph
$\Gamma(f)\subset X\times X$ so that the first projection $\pi_1\colon \Gamma(f) \to X$ is a
birational morphism. The locus where the inverse $\pi_1^{-1}$ is not well-defined is the
indeterminacy locus $\Ind(f) \subset X$: it is a closed subset of codimension at least $2$. Let
$\Delta$ be the diagonal of $X\times X$, we set the intersection schemes
$F_n(f) \coloneqq \Gamma(f^n)\times_{X\times X} \Delta$ whose support is the intersection of
$\Gamma(f^n)\cap\Delta$.

In general, the indeterminacy locus is not empty and the set of periodic points has infinite
cardinality, so it is natural to consider the set of isolated actual periodic points. Let $U_n$ be
the open locus on which the length-$n$ orbit under $f$ is well defined. More precisely, we can
define $U_n$ by induction,
$$
  U_1 =X\backslash \Ind(f) ,\qquad U_{n+1} = \left\{x\in U_n, f^n(x)\not\in \Ind(f) \right\}.
$$

The scheme of isolated periodic points of period $n$, denoted by $\IPer_n(f)$, is the open and
closed subscheme of \(\Per_n(f)\) supported on the union of its zero-dimensional irreducible
components contained in \(U_n\).

We introduce the ways to count points in $\IPer_n(f)(K)$: without multiplicity denoted by
$\NIPer_n(f)$; with multiplicity denoted by $\NIPerM_n(f)$. We refer to
Section~\ref{section-periodic-pts} for definitions and details. We are interested in the asymptotic
behavior of counts.

\begin{quest}
  What is the growth rate of $\NIPer_n(f)$ (resp. $\NIPerM_n(f)$) as $n\to \infty$.
\end{quest}

When $X$ is sufficiently regular, for instance when it is smooth, we may use intersection theory to
relate $\NIPer_n(f)$ to the
Lefschetz number of $f^n$, which is defined as the intersection number
\[
  L(f^n)\coloneqq\Gamma(f^n)\cdot \Delta \in \Z
  ~.
\]
The strategy is then to use excess formulas in order to write
$L(f^n)$ as a sum of suitably defined multiplicities assigned to each irreducible component of
$\Gamma(f^n)\cap \Delta$.

A computation due to Gromov allows one to control the intersection number $\Gamma(f^n)\cdot \Delta$
in terms
of the dynamical degrees of $f$, defined below. These numbers control the dynamical complexity
of an algebraic dynamical system and are strongly related to entropy when the base field is
metrized.

We may thus achieve our goal once we can show that
$\NIPer_n(f) \sim L(f^n)$. Proving such a comparison requires overcoming several difficulties.
\begin{enumerate}
  \item The multiplicities of the components of $\Gamma(f^n)\cap \Delta$ can be very high and grow
    uncontrollably with $n$.
  \item Positive-dimensional components of $\Gamma(f^n)\cap \Delta$ contribute to the computation of
    $L(f^n)$
    and, except in favorable situations (e.g. when $X$ is a homogeneous space), these contributions
    may be positive or negative.
  \item Some components of $\pi_1(\Gamma(f^n)\cap \Delta)$ may intersect the indeterminacy locus,
    and one also needs
    to bound their multiplicities.
\end{enumerate}

Our first result is a general upper bound on $\NIPer_n(f)$ in terms of degrees of $f^n$.
To state our theorem, fix an ample line bundle $L\to X$, and recall that for any
$p\in\{0,\ldots,m\}$ the degree
of $f$ with respect to $L$ is defined by
\[
  \deg_{L,p}(f)
  =
  \pi_1^\ast c_1(L)^{m-p}\cdot \pi_2^\ast c_1(L)^p\cdot \Gamma(f).
\]
A theorem of Dinh--Sibony in characteristic zero, and of
\cite{dang2020degrees} and \cite{truong2020relative} over any field,
states that the sequence $\deg_{L,p}(f^n)$ is submultiplicative up to a uniform constant, so that
one may define the dynamical degrees
\[
  \lambda_p(f)\coloneqq
  \lim_{n\to\infty}\deg_{L,p}(f^n)^{1/n}
\]
by Fekete's lemma. Recall that changing the polarization may change the degrees only up to a
multiplicative constant
which is independent of the iterate, so that the dynamical degrees do not depend on this choice.

The study of the number of perfect points, which generalize isolated periodic points, dates back to
Fulton. \cite{vu2022exotic} called the perfect but not periodic points are exotic periodic points
and investigated the growth of the related count for complex dynamics using density of currents. To
state the general version of our theorem, we first introduce the notion of a perfect point. Assume
that $X$ is smooth. Following Fulton, we consider the blow-up of the diagonal
$\pi\colon \widetilde{X\times X} \to X\times X$ and look at the schematic intersection of the
exceptional divisor $E\coloneqq\pi^{-1}(\Delta)$ with the strict transform $\widetilde{\Gamma(f^n)}$
of the graph. The crucial observation is that $E$ is a hypersurface, so the intersection
$E\cap \widetilde{\Gamma(f^n)}$ has pure dimension $m-1$. We call $P =(x,x) \in \Delta(K)$ a
\textit{perfect point} of period $n$ if the fiber $E\cap \widetilde{\Gamma(f^n)} \to \Delta$ over
$P$ has full dimension i.e. contains $E_P=\pi^{-1}(P)$. We emphasize that all isolated points are
perfect. Denote by $\Perfect_n(f)$ the number of perfect points of $\Gamma(f^n)$ with geometric
multiplicity of $E_P$ in $E\cap \widetilde{\Gamma(f^n)}$. We have
\[
  \NIPerM_n(f) \leq \Perfect_n(f).
\]
For details, we refer to Section~\ref{section-periodic-pts}.

\begin{theo} \label{theoremA}
  Let $X$ be an $m$-dimensional smooth projective variety over an algebraically
  closed field. Then there exists a constant $C>0$ such that for any dominant
  rational map $f\colon X\dasharrow X$ and any integer $n\in\N$, we have
  \begin{align}
    \Perfect_n(f)\le C\max_{0\le p\le m}\deg_{L,p}(f^n). \label{main-theoA}
  \end{align}
  In particular, we have
  \begin{align}
    \limsup_{n\to\infty}\frac{1}{n}\log \Perfect_n(f)
    \le \log\max_{0\le p\le m}\lambda_p(f). \label{cor-theoA}
  \end{align}
\end{theo}

Since
\[
  \NIPer_n(f)\leq \NIPerM_n(f)=N^{\text{geom}}_n(f) \leq \Perfect_n(f),
\]
we thus obtain an upper bound for all the quantities. A version for the count of isolated periodic
points with multiplicity of inequality~\eqref{cor-theoA},
\[
  \limsup_{n\to\infty}\frac{1}{n}\log \NIPerM_n(f)
  \le \log\max_{0\le p\le m}\lambda_p(f),
\]
was first proved by Dinh, Nguyen, and Truong \cite{dinh2017growth} for meromorphic self-maps of
compact K\"ahler manifolds, using the analytic theory of densities of positive closed currents. Vu
\cite{vu2022exotic}, \cite{vu2019densities} proved \eqref{main-theoA} using by bounding the mass of
strict transforms of currents under blow-ups using local analytic estimates. He went further and
used the theory of density of currents to obtain lower bounds for \(\Perfect_n(f)\) for special
classes of meromorphic maps (including those with dominant topological degree).

Our proof relies exclusively on Fulton's intersection theory and works over algebraically closed
fields of arbitrary characteristic.

The quantity $\log \max_{0\leq p\leq m}\lambda_p(f)$ is usually
referred to as the algebraic entropy of $f$; see~\cite{dinh2017growth,bellon1999algebraic}.
When the base field is metrized and complete, it is known to give an upper bound for
the topological entropy of the map $f$ induced on the analytification of $X$:
in the archimedean case, this is a theorem of Dinh and Sibony~ \cite{dinh2005entropy}, and in the
non-archimedean
case, it is due to Favre, Truong, and Xie \cite{favre2022topological}. Whether one can directly
relate the number of
isolated periodic points to the topological entropy is a delicate issue.

The key step of Theorem~\ref{theoremA} (Theorem~\ref{complexity-inequality} and
Proposition~\ref{proposition-the-three-equivalence-volume-estimates}) is to prove that the volume of
$E\cap \widetilde{\Gamma(f^n)}$ in $\widetilde{X\times X}$ is controlled uniformly by the volume of
$\Gamma(f^n)$ in $X\times X$.

\medskip

We now turn to the problem of obtaining lower bounds on $\NIPer_n(f)$, which is more
subtle because one must produce sufficiently many isolated periodic
points. Standard techniques rely heavily on pluripotential theory and only work
over the complex numbers.
Our second main result gives the exact asymptotic growth for birational maps of
surfaces in characteristic zero.

\begin{theo}\label{theoremB}
  Let $f\colon X\dasharrow X$ be a birational map of a normal projective surface
  defined over an algebraically closed field of characteristic $0$. Assume that $\lambda_1(f)>1$.
  Let $Q_n$ be one of $\NIPer_n(f)$, $\NIPerM_n(f)$. Then
  \begin{align*}
    Q_n(f)= \lambda_1(f)^n+o(\lambda_1(f)^n).
  \end{align*}
\end{theo}

We first prove the upper bound for the number of graph-periodic points
$$
  \ \limsup_{n\to +\infty}\frac{\NIPerM_n(f)}{\lambda_1(f)^n} \le 1
$$
(Theorem~\ref{upper-bound-surface}) by refining the arguments leading to Theorem~\ref{theoremA}.
Keeping the same notation as above, we need to find an upper bound for the sum of the multiplicities
of the components of $\widetilde{\Gamma(f^n)} \cap E$ that are horizontal with respect to the
fibration $\pi\colon E\to \Delta$. Note that any such component projects onto a curve of periodic
points in $X$. We obtain this bound by combining Cantat's finiteness theorem
\cite{cantat2010invariant} for periodic curves with an explicit calculation showing that the
multiplicities of these horizontal curves are bounded in terms of the iterates. Let us also note the
Saito's formula and works in \cite{iwasaki2010periodic}.

The proof of the lower bound (Theorem~\ref{theorem-main-lower-bound})
$$
  \liminf_{n\to \infty} \frac{\NIPer_n(f)}{\lambda_1(f)^n}\ge 1
$$
requires two independent steps. Suppose first that the base field is the field of complex numbers.
Bedford and Diller \cite{bedford2005energy} introduced an "energy" condition quantifying how fast
points of indeterminacy of $f^{-1}$ can approach $\Ind(f)$, and constructed an invariant measure
under this condition. It was later proved by Diller, Dujardin and Guedj \cite{ddg3dynamics} that
saddle periodic points are equidistributed with respect to this measure, and their arguments also
yield the required lower bounds. Jonsson and Reschke \cite{jonsson2018complex} proved that, if $f$
is defined over a number field, then every base change of $f$ to $\IC$ satisfies the Bedford--Diller
energy condition. Using the theory of Moriwaki heights, we extend their arguments to any finitely
generated field of characteristic zero, and show that we may find a specialization to $\IC$ such
that the induced map satisfies the Bedford--Diller energy condition.

Although there are examples of complex birational map which do not satisfy the Bedford-Diller energy
condition, see~ \cite{diller2009regularity} as well as \cite{favre1998points}. Jonsson and Reschke
\cite{jonsson2018complex} proved that, when \(f\) is defined over a number field, then every base
change of \(f\) to \(\IC\) does. Our argument builds on their result and uses the theory of Moriwaki
heights. We first choose a finitely generated field \(K'\) over which both \(X\) and \(f\) are
defined, and then prove that, for "most" (in a measurable sense) embeddings of \(K'\) in \(\IC\),
the induced complex map satisfies the Bedford–Diller energy condition.

\subsection*{Notes and references.}

The general viewpoint that dynamical zeta function should be encoded by the number of periodic
points \cite{artin1965periodic} shows the study of counts is valuable.

For dimension one, especially for maps of $\IP^1$, we refer to Silverman's book
\cite{silverman2010arithmetic}. Over fields of positive characteristic, purely inseparable phenomena
have to be separated from the separable dynamics.

For polarized endomorphisms, the picture is much more rigid. If $f^\ast L\simeq L^{\otimes q}$ with
$q>1$, then the degrees of the iterates are controlled by the single scaling factor $q$. Fakhruddin
proved, in particular, that periodic points are Zariski dense for polarizable self-maps of
projective varieties over algebraically closed fields \cite{fakhruddin2003questions}. Over $\IC$,
the basic model is the case of holomorphic endomorphisms of projective space, where Briend--Duval
proved equidistribution results for periodic points \cite{briend1999periodic}.

For general dominant rational or meromorphic self-maps, the main difficulty
is that fixed loci of iterates may have positive-dimensional components. In the
compact K\"ahler setting, Dinh--Nguyen--Truong proved that the exponential
growth rate of isolated periodic points, counted with multiplicity, is bounded
above by the algebraic entropy \cite{dinh2017growth}. Their proof uses the
theory of densities of positive closed currents developed by Dinh--Sibony
\cite{dinh2018density}. Perfect, or exotic, periodic points in this general
intersection-theoretic setting were studied by Vu \cite{vu2022exotic}.

When the topological degree dominates the other dynamical degrees, stronger
results are available. For meromorphic self-maps of compact K\"ahler manifolds
with dominant topological degree, Dinh--Nguyen--Truong proved
equidistribution results for repelling periodic points and asymptotic formulas
in terms of the topological degree \cite{dinh2015equidistribution}. This case
is closer in spirit to the classical expanding situation.

For birational maps of surfaces, the first dynamical degree is the main source
of exponential growth. For polynomial diffeomorphisms of $\IC^2$,
Bedford--Lyubich--Smillie proved equidistribution results for periodic points
\cite{bedford1993distribution}. For birational maps of $\IP^2$, Favre studied
the growth of non-critical periodic points \cite{favre1998points}. The
dynamics of bimeromorphic surface maps was developed systematically by
Diller--Favre \cite{diller2001dynamics}. Xie proved that, for birational maps
of smooth projective surfaces with first dynamical degree $>1$, non-critical
periodic points are Zariski dense \cite{xie2012periodic}; see also his thesis
\cite{xie2014dynamique}. Periodic curves play a special role in this setting:
Cantat proved finiteness results for invariant hypersurfaces
\cite{cantat2010invariant}, and the comparison between Lefschetz numbers and
isolated periodic points in the bimeromorphic K\"ahler case can be expressed
using the Saito--Iwasaki--Uehara formula; see
\cite{dinh2017growth,iwasaki2010periodic}. For surface maps of small
topological degree, Diller--Dujardin--Guedj constructed invariant measures and
proved equidistribution of saddle periodic points under an energy condition
\cite{ddg3dynamics}. The problem of equidistribution of saddle periodic points are also studied in several cases by \cite{luo2026equidistribution},\cite{luo2026measures} and \cite{luo2025equidistribution}.
\bigskip

\noindent{\bf Plan of the paper.}
The article is organized as follows. In Section~\ref{section2}, we recall the
intersection-theoretic tools and develop the volume inequalities. In
Section~\ref{section-periodic-pts}, we present the settings of counts problems.
Section~\ref{section3} proves Theorem~\ref{theoremA} and the upper bound for birational surface
maps. In Section~\ref{section4}, we recall the relevant parts of Moriwaki height theory, establish
the needed height estimates, and prove the lower bound in Theorem~\ref{theoremB}.

\medskip

\noindent{\bf Acknowledgments.} I would like to thank my advisor Charles Favre
for his supervision and for his help with the writing of this paper. I would
also like to thank the CMLS Laboratory for its support. I would like to thank China Scholarship
Council (CSC) as well as University of Science and Technology of China (USTC) for financial support.

\medskip

\noindent{\bf AI disclosure.}
Chat-GPT Sol 6.1 has been used to produce a proofreading report. The ideas in the appendix were also
suggested by Chat-GPT.Astra. The project is part of the Phd of the author and started in 2025. All
results in the introduction have been obtained without any AI-assistance. The entire manuscript is
human-written.
\section{Intersection theory and volumes}
\label{section2}

\subsection{Chow groups and Fulton's intersection theory}

We assume that varieties are integral and defined over an
algebraically closed field $K$. We use the conventions of Fulton's
intersection theory \cite[Chapters~1--2]{fulton2013intersection}.

\subsubsection{Chow groups}
Let $X$ be a normal projective variety of dimension $m$. For any integer $k$, the
group $Z_k(X)$ of $k$-cycles is the free abelian group generated by the
$k$-dimensional subvarieties of $X$. Let $W\subset X$ be an integral
$(k+1)$-dimensional subvariety and let $r\in K(W)^\times$. For every prime
divisor $D\subset W$, write $r=a/b$ in the one-dimensional local domain
$\cO_{W,\eta_D}$ and set
\[
  \ord_D(r)
  \coloneqq
  \ell_{\cO_{W,\eta_D}}\bigl(\cO_{W,\eta_D}/(a)\bigr)
  -\ell_{\cO_{W,\eta_D}}\bigl(\cO_{W,\eta_D}/(b)\bigr).
\]
This integer is independent of the representation $r=a/b$. The principal
divisor of $r$ on $W$ is the finite sum
\[
  \div_W(r)
  \coloneqq
  \sum_{\substack{D\subset W\\ \codim_W(D)=1}}
  \ord_D(r)[D]\in Z_k(W),
\]
which we also regard as a cycle on $X$. Two cycles $\alpha,\beta\in Z_k(X)$
are \emph{rationally equivalent}, written $\alpha\sim\beta$, if there exist
finitely many $(k+1)$-dimensional subvarieties $W_i\subset X$ and rational
functions $r_i\in K(W_i)^\times$ such that
\[
  \alpha-\beta=\sum_i\div_{W_i}(r_i).
\]
The
Chow group $A_k(X)$ is the quotient of $Z_k(X)$ by rational equivalence.
If $Z\subset X$ is a subvariety of dimension $k$, we denote by $[Z]$ its class
in $A_k(X)$.

Since $X$ is integral,
the map $n \mapsto n [X]$ induces an isomorphism
$A_m(X)\simeq \Z$. The group $A_{m-1}(X)$ is the group of Weil divisors modulo
rational equivalence. When $X$ is normal, it identifies with the divisor class
group $\operatorname{Cl}(X)$. Any line bundle $L\to X$
defines a class $c_1(L)\in A_{m-1}(X)$. When $X$ is smooth, this agrees with
the usual identification $A_{m-1}(X)\simeq \Pic(X)$.

Any $0$-cycle $\alpha$ can be written in the form
$\alpha = \sum a_i [p_i]$, where the $p_i$ are closed points. We set
$\deg(\alpha)\coloneqq \sum_i a_i$.
We obtain a morphism $\deg \colon Z_0(X) \to \Z$
which descends to a morphism
$\deg \colon A_0(X) \to \Z$.

Let $Z$ be a pure $k$-dimensional closed subscheme of $X$, and let
$Z_1,\cdots,Z_t$ be all the irreducible components of the underlying topological space $|Z|$. The
\textit{geometric multiplicity} of $Z_i$ in $Z$ is the
length of the Artinian local ring $\cO_{Z,\eta_i}$:
\[
  \ell(Z_i,Z)\coloneqq\ell_{\cO_{Z,\eta_i}}(\cO_{Z,\eta_i}), \label{geometric-multiplicity}
\]
where $\eta_i$ is the generic point of $Z_i$.

Let $K(Z_i)$ be its residue field. Since the maximal ideal $\mathfrak m_i$ of
$\cO_{Z,\eta_i}$ is nilpotent, additivity of length gives
\[
  \ell(Z_i,Z)=\sum_{s\geq 0} \dim_{K(Z_i)} \left( \mathfrak{m}_i^s/\mathfrak{m}_i^{s+1}\right),
\]
where the sum is finite; see
\cite[Tags~00IV and~00IY]{stacks-project}.
The fundamental cycle of $Z$ is then
\[
  [Z]=\sum_{i=1}^t \ell(Z_i,Z)[Z_i]\in A_k(X).
\]

\subsubsection{Proper push-forward and flat pull-back}
Let $X$ and $Y$ be projective varieties, and let $f\colon X\to Y$ be a
morphism. Since $f$ is proper, for any irreducible subvariety $V\subset X$,
the image $f(V)$ is a subvariety of $Y$. We set
\[
  f_\ast[V]=
  \begin{cases}
    0 & \text{if } \dim f(V)<\dim V,\\
    [K(V):K(f(V))]\,[f(V)] & \text{if } \dim f(V)=\dim V.
  \end{cases}
\]
This defines a push-forward morphism on cycles which descends to the Chow groups
\cite[Chapter~1]{fulton2013intersection},
\[
  f_\ast\colon A_k(X)\longrightarrow A_k(Y).
\]

Assume now that $f\colon X\to Y$ is flat of relative dimension $e$. For any
subvariety $V\subset Y$, the scheme-theoretic inverse image $f^{-1}(V)$ has
dimension $\dim V+e$, and we set
\[
  f^\ast[V]\coloneqq [f^{-1}(V)].
\]
This extends by linearity and descends to a morphism
\[
  f^\ast\colon A_k(Y)\longrightarrow A_{k+e}(X).
\]
We refer to \cite[Chapter~1]{fulton2013intersection} for details.

\subsubsection{Intersection with divisors}

Let $X$ be a normal projective variety of dimension $m$, and let $D$ be a Cartier divisor on $X$.
Let
$W\subset X$ be a $k$-dimensional subvariety, and let $j\colon W\hookrightarrow X$
be the inclusion. We define
\[
  D\cdot W \coloneqq
  j_\ast\bigl[\div_W(s)\bigr]\in A_{k-1}(X),
\]
where $s$ is any non-zero rational section of $\cO_X(D)|_W$. This class is
independent of the choice of section $s$. The morphism
$D\cdot -\colon Z_k(X)\to A_{k-1}(X)$ descends to a morphism of Chow groups. One can show that this
morphism depends only on the
linear equivalence class of $D$, so we may write
\[
  [D]\cdot -\colon A_k(X)\longrightarrow A_{k-1}(X),
\]
where $[D]\in A_{m-1}(X)$ corresponds to the linear equivalence class of $D$.

If $f\colon X'\to X$ is a morphism between projective varieties, then for any cycle
$\alpha\in A_\ast(X')$, the projection formula holds
\[
  f_\ast(f^\ast D\cdot \alpha)=D\cdot f_\ast\alpha.
\]
If $f$ is flat, then
\[
  f^\ast(D\cdot \beta)=f^\ast D\cdot f^\ast\beta
\]
for every $\beta\in A_\ast(X)$. Both identities follow from the compatibility
of Cartier divisors with proper push-forward and flat pull-back
\cite[Section~2.3]{fulton2013intersection}.

\subsubsection{Vector bundles and Chern classes}
We recall the facts about projective bundles that will be used below; see
\cite[Chapter~3]{fulton2013intersection} and
\cite[Section~9.3]{eisenbud3264}.
Let $N\to X$ be a vector bundle of rank $r$, and denote by $E=\IP(N)$ the
projective bundle of lines in $N$. Set $e = r-1$, the dimension of a fiber of $E$. Write
$p\colon E\to X$ for the projection,
and $E_x=p^{-1}(x)$ for the fiber over $x$. We use the convention that
$\cO_E(-1)$ is the tautological line subbundle of $p^\ast N$. Its dual line bundle
$\mathcal{O}_E(1)$ is relatively ample over $X$. In particular,
if $C\subset E_x$ is an integral curve contained in a fiber of
$p\colon E\to X$, then
\[
  h\cdot C
  =
  \deg\bigl(\mathcal{O}_E(1)|_C\bigr)
  >0,
  \qquad
  h=c_1\bigl(\mathcal{O}_E(1)\bigr).
\]
Note that $p$ is proper and flat; hence the push-forward and the pull-back are well defined.

\begin{lemma}\label{lemma-relative-ample}
  For any ample line bundle $L$ on $X$, there exists $N\in \IN$ large enough such that for $n\geq N$
  the line bundle $\cO_E(1) + n p^\ast L$ is ample.
\end{lemma}
\begin{proof}
  The argument is the projective-bundle case of
  \cite[Example~1.7.2 and Proposition~1.7.10]{lazarsfeld2017positivityI}.

  Since \(L\) is ample, there is an integer $n$ and a surjection
  \[
    \mathcal O_X^{\oplus s}\twoheadrightarrow
    N^\vee\otimes L^{\otimes n}.
  \]
  Dualizing gives an inclusion of vector bundles
  \[
    N\otimes L^{\otimes(-n)}
    \hookrightarrow \mathcal O_X^{\oplus s}.
  \]
  It induces a closed immersion of projective bundles of lines
  \[
    j\colon E=\IP(N)
    \cong \IP(N\otimes L^{\otimes(-n)})
    \hookrightarrow \IP^{s-1}\times X.
  \]
  Under this identification, the tautological line on
  \(\mathbb P(N\otimes L^{\otimes(-n)})\) is
  \(\mathcal O_E(-1)\otimes p^*L^{\otimes(-n)}\). Hence
  \[
    j^*\operatorname{pr}_1^*\mathcal O_{\IP^{s-1}}(1)
    \cong \mathcal O_E(1)\otimes p^*L^{\otimes n}.
  \]

  The exterior product
  \[
    \operatorname{pr}_1^*\mathcal O_{\mathbb P^{s-1}}(1)
    \otimes \operatorname{pr}_2^*L
  \]
  is ample on \(\mathbb P^{s-1}\times X\). Its restriction via \(j\)
  is \(\mathcal O_E(1)\otimes p^*L^{\otimes n+1}\), so this line bundle
  is ample. Thus \(N=n+1\) suffices.
\end{proof}
\begin{lemma}\label{lemmaprojectionformulaforprojectivebundle}
  For any class $\alpha\in A_\ast(X)$, we have
  \[
    p_\ast(p^\ast\alpha\cdot h^j)=
    \begin{cases}
      0 & \text{if } j<e,\\
      \alpha & \text{if } j=e .
    \end{cases}
  \]
\end{lemma}

\begin{proof}
  By the projection formula, we have
  \[
    p_\ast(p^\ast\alpha\cdot h^j)=\alpha\cdot p_\ast(h^j),
  \]
  so it is enough to prove $p_\ast(h^j)=0$ when $j<e$ and
  $p_\ast(h^{e})=[X]\in A_m(X)$.

  If $j<e$, then $h^j$ is represented by a cycle of dimension
  $\dim X+e-j>\dim X$, while its image under $p$ is contained in $X$.
  By the definition of proper push-forward, we have $p_\ast(h^j) =0$.

  When $j=e$, write $p_\ast(h^{e})=a[X]$ for some $a\in \Z$. It suffices
  to check on a closed point $x\in X$:
  \[
    a[x]=p_\ast\left(h^{e}\cdot[p^{-1}(x)]\right)
    =p_\ast\left(c_1(\cO_{\IP^{e}}(1))^{e}\right)=[x].
  \]
  This gives $a=1$.
\end{proof}

We recall the following decomposition theorem for the Chow groups of a projective bundle.

\begin{theorem}[{\cite[Theorem~3.3]{fulton2013intersection}}]\label{structure-of-chow-ring-of-projective-bundle}
  For each $k\in\N$, the flat pull-back morphism $p^\ast\colon A_{k}(X)\to A_{k+e}(E)$ is injective
  and there is a canonical isomorphism
  \[
    A_k(E) \simeq \bigoplus_{i=0}^{e}A_{k-e+i}(X).
  \]
  More precisely, each element $\beta \in A_k(E)$ can be uniquely expressed in the form
  \[
    \beta = \sum_{i=0}^e p^\ast(\alpha_i) h^i,
  \]
  where $\alpha_i\in A_{k-e+i}(X)$.
\end{theorem}

\begin{definition}
  The Segre classes of $N$ are defined by
  \[
    s_i(N)=p_\ast(h^{e+i})\in A_{m-i}(X),\qquad i\geq 0.
  \]
  The total Segre class is
  \[
    s(N)=\sum_{i\geq 0}s_i(N)\in A_\ast(X) .
  \]
  Note that the intersections with Segre classes are well defined. For any $\alpha\in A_\ast(X)$, we
  set
  \[
    s_i(N)\cdot \alpha=p_\ast(h^{e+i}\cdot p^\ast\alpha).
  \]

  We then define the Chern classes by the identity
  \[
    c_t(N)s_t(N)=1,
  \]
  where $s_t(N)=\sum_{i\geq 0}s_i(N)t^i$ and
  $c_t(N)=\sum_{i \geq 0} c_i(N)t^i$ are formal power series. Thus $c_0(N)=1$,
  $c_1(N)=-s_1(N)$, and, for $n\geq 1$,
  \[
    c_n(N)=-\sum_{i=1}^{n}s_i(N)c_{n-i}(N).
  \]
  In particular, we have $c_i(N)=0$ for $i>r=e+1$ and $s_i(N) = 0$ for $i >m$.
\end{definition}

Since the Chern classes are polynomials in the Segre classes, their intersections are also well
defined. In fact, they can be described as follows.
For an irreducible subvariety $Z$ of $X$, we have
$c_i(N) \cdot [Z] = j_\ast c_i(N|_Z)$,
where $j:Z\hookrightarrow X$ is the inclusion map.

The Chern classes are compatible with flat pull-back:
\[
  c_i(f^\ast N)\cdot f^\ast\alpha=f^\ast(c_i(N)\cdot\alpha)
\]
for any flat morphism $f\colon X'\to X$. They also satisfy the Whitney formula:
for every exact sequence of vector bundles
\[
  0\to N'\to N\to N''\to 0
\]
one has the following relation of formal polynomials
$$
  c_t(N)=c_t(N')c_t(N'').
$$

\subsubsection{Normal cones and Segre classes}
\label{subsubsection-segre-normal-cones}

Cones generalize vector bundles. Although a general cone has no Chern
classes in the usual vector-bundle sense, its projective completion carries
a tautological line bundle used to define its Segre class; see
\cite[Sections~4.1--4.3]{fulton2013intersection}.

Let $X$ be a normal projective variety of dimension $m$. If $Z\subset X$ is a closed
subscheme defined by an ideal sheaf $\mathcal I$, then
\[
  C_ZX\coloneqq
  \operatorname{Spec}_Z\left(\bigoplus_{n\geq 0}\mathcal I^n/\mathcal I^{n+1}\right)
\]
is the normal cone of $Z$ in $X$. It is a cone of pure dimension $m$ over
$Z$; see \cite[Appendix~B.6.1 and~B.6.6]{fulton2013intersection}.

For a cone $C=\operatorname{Spec}_Z(\mathcal S^\bullet)$ over $Z$, let
\[
  \overline C=\IP(C\oplus 1)
  \coloneqq\operatorname{Proj}_Z(\mathcal S^\bullet[u]),
  \qquad q\colon\overline C\longrightarrow Z,
\]
be its projective completion, where $u$ has degree $1$ and $1$ denotes the
trivial line bundle. Put $\overline h=c_1(\cO_{\overline C}(1))$. Following
\cite[Sections~4.1--4.2]{fulton2013intersection}, we define
\[
  s(C)\coloneqq q_\ast\left(\sum_{i\geq 0}\overline h^{\,i}\cdot
  [\overline C]\right),
  \qquad s(Z,X)\coloneqq s(C_ZX)\in A_\ast(Z).
\]
The sum is finite by dimension. We write $\{s(Z,X)\}_k$ for its component in
$A_k(Z)$. If $Z$ has pure dimension $z$, we set
$s_i(Z,X)\coloneqq\{s(Z,X)\}_{z-i}$, so that $i$ measures codimension on
$Z$, as for the vector-bundle Segre classes above. In particular,
$C_XX$ is the zero cone and $s(X,X)=[X]$.

If $Z\neq X$, the blow-up of $X$ along $Z$ is
\[
  \pi\colon\widetilde X
  =\operatorname{Proj}_X\left(\bigoplus_{n\geq 0}\mathcal I^n\right)
  \longrightarrow X.
\]
Its exceptional divisor is $E=\IP(C_ZX)$; see
\cite[Appendix~B.6.3]{fulton2013intersection}. Denote its projection by
$p\colon E\to Z$ and put $h=c_1(\cO_E(1))$. Then
$\cO_E(1)=\cO_{\widetilde X}(-E)|_E$, and
\cite[Corollary~4.2.2]{fulton2013intersection} gives the equivalent formula
\[
  s(Z,X)=p_\ast\left(\sum_{i\geq 0}h^i\cdot[E]\right).
\]

Recall that a closed embedding $Z\hookrightarrow X$, defined by an ideal
sheaf $\mathcal I$, is a regular embedding of codimension $d$ if $\mathcal I$
is locally generated by a regular sequence $f_1,\ldots,f_d$; that is, each
$f_j$ is a non-zero-divisor modulo $(f_1,\ldots,f_{j-1})$; see
\cite[Appendix~B.7]{fulton2013intersection}. In this case,
$\mathcal I/\mathcal I^2$ is locally free of rank $d$, and the normal cone
$C_ZX$ is the normal bundle $N_{Z/X}$. Consequently,
\[
  s(Z,X)=c(N_{Z/X})^{-1}\cdot[Z],
\]
which agrees with the Segre class of the vector bundle $N_{Z/X}$ by
\cite[Proposition~4.1(a)]{fulton2013intersection}.

\subsubsection{Multiplicity along a subvariety}\label{section-Multiplicity along subvariety}

For an irreducible subvariety $V\subsetneq X$ of codimension $d$, the
multiplicity of $X$ along $V$, denoted by $e_V(X)$, is defined by
\[
  \{s(V,X)\}_{m-d}=e_V(X)[V]
  =p_\ast\left(h^{d-1}\cdot[\IP(C_VX)]\right),
\]
where $p$ and $h$ now refer to $\IP(C_VX)$. Equivalently, if
$A=\cO_{X,\eta_V}$ and $\mathfrak m_A$ is its maximal ideal, then
$e_V(X)=e_{\mathfrak m_A}(A)$ is the Hilbert--Samuel multiplicity defined by
\[
  \ell_A(A/\mathfrak m_A^t)
  =\frac{e_{\mathfrak m_A}(A)}{d!}t^d+O(t^{d-1})\qquad(t\to\infty).
\]
See \cite[Section~4.3 and Example~4.3.1]{fulton2013intersection}.

More generally, if $V$ is an irreducible component of a proper closed
subscheme $Z\subset X$, the multiplicity of $X$ along $Z$ at $V$ is the coefficient of $[V]$ in
$s(Z,X)$ near $\eta_V$, denoted by $(e_Z X)_V$; see
\cite[Example~4.3.4]{fulton2013intersection}. This should be distinguished from the geometric multiplicity $\ell(V,Z)$ defined in Section~\ref{geometric-multiplicity}. Lech’s inequality \cite{ma2023uniform} gives $(e_Z X)_V \leq m!\, \ell(V,Z) e_V(X)$.

\subsubsection{Intersection products}

Let $X$ be a normal projective variety of dimension $m$. Let
$i\colon V\hookrightarrow X$ be a regular embedding of codimension $d$, with
normal bundle $N_{V/X}$. Let $W$ be a normal projective variety of dimension $n$, and
let $f\colon W\to X$ be a morphism. Form the fiber square
\[
  \begin{tikzcd}
    Z=V\times_X W \arrow[r, "j", hook] \arrow[d, "g"']
    & W \arrow[d, "f"] \\
    V \arrow[r, "i"', hook] & X .
  \end{tikzcd} \label{intersection-square}\tag{$\star$}
\]

One can check that the normal cone $C_ZW$ is naturally a closed
subcone of the rank $d$ vector bundle $N\coloneqq g^\ast N_{V/X}$ . Fulton defines the intersection
product of $W$ with $V$ over $X$
by
\[
  V\cdot_X W
  =\left\{c(N)\cdot s(Z,W)\right\}_{n-d}\in A_{n-d}(Z),
\]
where $\left\{-\right\}_{n-d}: A_\ast(Z)=\bigoplus_k A_k(Z) \to A_{n-d}(Z)$ is the projection to the
$n-d$-th factor. See \cite[Sections~6.1--6.2]{fulton2013intersection}. Note that the product class
depends only on the rational equivalence class $[W]$.

More generally, Fulton defines a refined Gysin homomorphism
\[
  i^!\colon A_k(W)\longrightarrow A_{k-d}(Z),
\]
such that
$$
  V\cdot_X W=i^![W].
$$

\medskip

Let $C_1,\cdots, C_r$ be the irreducible components of $C_ZW$ with geometric multiplicities
$l_1,\cdots, l_r$. Then the Segre class is $s(Z,W) = \sum l_i s(C_i)$. One can write
\[
  V\cdot_X W
  =\sum l_i \left\{c(N)\cdot s(C_i)\right\}_{n-d}.
\]

When $W$ is a subvariety of $X$, we simply denote by $V\cdot W$ the intersection product class
$V\cdot_X W$ in $A_{n-d}(Z)$.

Suppose that $V=D$ is an effective Cartier divisor. Then the
intersection product agrees with the intersection
with divisors introduced above.

\subsubsection{Proper intersections} \label{subsubsection-proper-intersection}
We use the notations as \eqref{intersection-square}. An irreducible component $C$ of the
intersection $Z = V\cap W$ is called proper if
$\dim C=n-d$. If every component is proper, then
\[
  V\cdot W
  =\sum_C i(C,V\cdot W;X)[C]\in A_{n-d}(Z),
\]
where the sum runs over the irreducible components of $Z$ and the intersection multiplicity
$$
  i(C,V\cdot W;X) = (e_Z W)_C
$$
is the multiplicity of $W$ along $Z$ at $C$; see \cite[Section~7.1]{fulton2013intersection}.

The intersection multiplicity can be expressed as Serre's intersection multiplicity formula. In
general, we have
\[
  1\leq i(C,V\cdot W;X)\leq \ell(C,Z)
\]
where $\ell(C,Z)$ is the geometric multiplicity.

The intersection multiplicity coincides with the geometric multiplicity
$i(C,V\cdot W;X) = \ell(C,Z)$ if $W$ is Cohen--Macaulay at the generic point of $C$; see
\cite[Proposition~7.1]{fulton2013intersection}.

Let $A= \cO_{Z,\eta_C}$ and $J$ the ideal in $A$ generated by ideal of $W$ in $X$. Then
$i(C,V\cdot W;X)=1$ if and only if $\ell(C,Z) =1$ if and only if $A$ is regular and
$J=\mathfrak m_A$. When $X$ is smooth, this is equivalent to $C$ being a transverse component. Here,
a proper component $C$ is called {transverse} if, at a general point $x\in C$, both $V$ and $W$ are
smooth and
\[
  T_xV+T_xW=T_xX.
\]

\subsubsection{Chow ring structure}

Suppose that $X$ is a smooth projective variety of dimension $m$. Set $A^d(X) = A_{m-d}(X)$ for each
$d\in \{0,\cdots,m\}$. Then the diagonal embedding
$\Delta\colon X\hookrightarrow X\times X$ is regular of codimension $m$ with normal bundle
$N_{\Delta/X\times X} \simeq T_X$. Let $p_1,p_2\colon X\times X\to X$ be the two projections. The
exterior product
is the homomorphism
\[
  \boxtimes\colon A^d(X)\otimes_\Z A^e(X)
  \longrightarrow A^{d+e}(X\times X)
\]
defined on classes of subvarieties by
\[
  [Z]\boxtimes[W]
  \coloneqq [Z\times W]
  =p_1^*[Z]\cdot p_2^*[W].
\]
The intersection product is the composition
\[
  A^d(X)\otimes_\Z A^e(X)
  \xrightarrow{\ \boxtimes\ } A^{d+e}(X\times X)
  \xrightarrow{\Delta^{!}} A^{d+e}(X),
\]
and it makes $A^\ast(X)=\bigoplus_d A^d(X)$ into a graded commutative ring.
The unit is the fundamental class
$[X]\in A^0(X)$. This product agrees with the intersection with Cartier
divisors defined above and with the intersection product whenever the
cycles meet properly.

The Chow ring of a projective bundle may be viewed as an extension of the Chow ring of $X$; see
\cite[Theorem~9.6]{eisenbud3264}. There is a refinement of
Theorem~\ref{structure-of-chow-ring-of-projective-bundle}.

\begin{theorem}[Projective bundle formula]
  \
  Let $X$ be smooth and let $N\to X$ be a vector bundle of rank $r$. For
  $E=\IP(N)$, there is a canonical isomorphism of graded rings
  \[
    A^\ast(E)\simeq A^\ast(X)[h]/(h^r+c_1h^{e}+\cdots+c_r),
  \]
  where $c_i=c_i(N)\in A^i(X)$ is the $i$-th Chern class of $N$.
\end{theorem}

\subsubsection{Pull-back operations on Chow groups}
Let $f\colon X\to Y$ be a morphism between smooth projective varieties of
dimensions $m$ and $n$. There is a pull-back of graded rings
\[
  f^\ast\colon A^\ast(Y)\longrightarrow A^\ast(X).
\]
To construct it, let
$\Gamma_f\colon X\hookrightarrow X\times Y$ be the graph embedding and let
$p_2\colon X\times Y\to Y$ be the second projection. Since $\Gamma_f$ is a
regular embedding of codimension $n$, one can define
\[
  f^*\alpha
  \coloneqq
  \Gamma_f^!\bigl(p_2^*\alpha\bigr),
  \qquad \alpha\in A^*(Y),
\]
where $p_2^*$ is the flat pull-back and $\Gamma_f^!$ is the refined Gysin
homomorphism. If $f$ is flat, this
pull-back agrees with the flat pull-back defined above and, for every Cartier
divisor $D$ on $Y$, it satisfies
\[
  f^\ast c_1(\cO_Y(D))=c_1(\cO_X(f^\ast D)).
\]
Thus the construction is compatible with both pull-back operations introduced
earlier. We refer to
\cite[Chapter~8]{fulton2013intersection} and
\cite[Theorem~1.23]{eisenbud3264} for the construction by refined Gysin maps.

These operations satisfy the projection formula. More precisely, if
$\alpha\in A^d(Y)$ and $\beta\in A^e(X)$, then
\[
  f_\ast\bigl(f^\ast\alpha\cdot \beta\bigr)
  =
  \alpha\cdot f_\ast\beta
  \in A^{d+e+n-m}(Y).
\]

\subsection{Numerical classes and positivity}
Let $X$ be a smooth projective variety of dimension $m$.

\subsubsection{Numerical classes}
We are mainly interested in intersection numbers, so it is natural to
work in the space of cycles modulo numerical equivalence. A class $\alpha\in A^d(X)$
is numerically trivial if $[Z]\cdot\alpha=0$ for every $d$-dimensional subvariety
$Z\subset X$. In this setting, we write
$Z^d(X)=Z_{m-d}(X)$. Let $N^d(X)$ denote the quotient
of $A^d(X)$ by numerical equivalence and set
\[
  Z^d_\R(X)\coloneqq Z^d(X)\otimes_\Z\R,
  \qquad
  N^d_\R(X)\coloneqq N^d(X)\otimes_\Z\R .
\]
The spaces $N^d_\R(X)$ are finite-dimensional, and
$N^0(X) \simeq A^0(X) \simeq \Z$; see
\cite[Theorem~2.5.1]{dang2020degrees} and \cite[Theorem~23.6]{milne2012lectures}.
The degree map induces an isomorphism $N^m(X)\simeq \Z$.
The intersection product on $A^\ast(X)$ induces a non-degenerate pairing
$N^d_\R(X)\times N^{m-d}_\R(X)\longrightarrow \R.$

\subsubsection{Positive classes}
Let us discuss some positivity notions on cycles, following
\cite{lazarsfeld2017positivityI,fulger2014positive}.

A cycle $\alpha = \sum_i a_i [Z_i]$ is effective if $a_i \ge 0$ for all $i$. A class
$\alpha\in N^d_\R(X)$ is effective if it is represented by an effective cycle.
The convex cone generated by all effective classes is denoted by $\Eff^d(X)$, and
its closure $\Psef^d(X)$
is the pseudo-effective cone. It has non-empty interior and is salient, i.e.,
$\pm \alpha \in \Psef^d(X)$ implies $\alpha=0$.
A class is big if it lies in the interior of $\Psef^d(X)$. For any
$\alpha,\beta\in N^d_\R(X)$, we write $\alpha\geq\beta$ if
$\alpha-\beta\in\Psef^d(X)$.

By duality, we say that a class $\alpha\in N^d_\R(X)$ is
nef if $\alpha\cdot\beta\geq 0$
for every pseudo-effective class $\beta\in N^{m-d}_\R(X)$.
The set of nef classes $\Nef^d(X)$ is a closed salient cone
dual to $\Psef^{m-d}(X)$.
When $d=1$, we have $\Nef^1(X)\subset \Psef^1(X)$, and the interior
of $\Nef^1(X)$ consists of ample classes, i.e., it is the positive cone generated
by all classes $c_1(L)$ with $L$ ample.

These cones behave well under the functorial operations defined above. Proper
push-forward sends effective cycles to effective cycles, hence sends
pseudo-effective classes to pseudo-effective classes. Flat pull-back also sends
effective cycles to effective cycles, hence preserves pseudo-effectivity. When
a pull-back morphism is available and satisfies the projection formula, for
any $\alpha\in \Nef^d(Y)$, the class
$f^\ast\alpha$ is nef on $X$: indeed, for every pseudo-effective class
$\beta\in N^{m-d}_\R(X)$, the projection formula gives
\[
  f^\ast\alpha\cdot\beta = \alpha\cdot f_\ast\beta \geq 0,
\]
because $f_\ast\beta$ is pseudo-effective.

\subsubsection{Volumes}\label{section--volume}
Let us fix an ample line bundle $L\to X$.
Note the function
$$
  \Psef^d(X) \to \IR_{\geq 0},\quad \alpha \mapsto \alpha\cdot c_1(L)^{m-d},
$$
is
the restriction to $\Psef^d(X)$ of a norm on $N^d_\R(X)$, see \cite[Section~2.3]{dang2021spectral},
\cite[Theorem~1.4]{fulgerLehmann2017dual}.
It follows that, given any norm $\|\cdot\|$ on $N^d_\R(X)$, there
exists a positive constant $C>1$ such that
\begin{align}
  C^{-1}\|\alpha\|\leq \alpha \cdot c_1(L)^{m-d} \leq C\|\alpha\|\label{norm-and-volume-comparison}
\end{align}
for any $\alpha \in \Psef^d(X)$.
To simplify notation, for any $\alpha\in \Psef^d(X)$ we shall write
\begin{align*}
  \Vol_L(\alpha)\coloneqq \alpha\cdot c_1(L)^{m-d} \geq 0.
\end{align*}
We call this number the volume of $\alpha$ with respect to $L$.

When $Z\subset X$ is a subscheme, we denote by $\Vol_L(Z)$ the volume of the fundamental cycle of
$Z$. If $Z$ has pure codimension $d$, then $\Vol_L(Z) =[Z]\cdot c_1(L)^{m-d}$.

In the complex projective setting, when $c_1(L)$ is represented by a K\"ahler
form $\omega$, for a subvariety $Z$, Wirtinger's theorem identifies $\frac{1}{m!} \Vol_L(Z)$ with
the volume of
the smooth locus of $Z$ for the canonical volume form induced by the Riemannian
metric $\omega|_Z$; see \cite[Chapter~III, Section~1.D]{demailly1997complex}.

\begin{lemma}\label{lemma-effective-and-norms}
  There exists a constant $C>0$, depending only on $(X,L,\|\cdot\|)$, such that if
  $\beta\in\Psef^d(X)$ and $\alpha\in N^d_\R(X)$ satisfy
  \[
    -\beta\leq\alpha\leq\beta,
  \]
  then
  \[
    \|\alpha\|\leq C\|\beta\|.
  \]
\end{lemma}
\begin{proof}
  Let $C_0>1$ be the constant in~\eqref{norm-and-volume-comparison}. The
  hypothesis says that $\beta+\alpha$ and $\beta-\alpha$ are pseudo-effective.
  Consequently,
  \begin{align*}
    \|\alpha\|
    & \leq \frac12\bigl(\|\beta+\alpha\|+\|\beta-\alpha\|\bigr)\\
    & \leq \frac{C_0}{2}
    \bigl(\Vol_L(\beta+\alpha)+\Vol_L(\beta-\alpha)\bigr)\\
    & =C_0\Vol_L(\beta)
    \leq C_0^2\|\beta\|.
  \end{align*}
  Thus one may take $C=C_0^2$.
\end{proof}

\subsection{Volume inequalities}

Let $X$ be a smooth projective variety of dimension $m$, and let $L$ be an
ample line bundle on $X$. Let $V$ be a smooth irreducible subvariety of codimension $e$, then
$i\colon V\hookrightarrow X$ is a regular embedding. Let $\pi\colon \widetilde X\to X$ be the
blow-up of $X$ along $V$. The exceptional divisor $E$ is the projective normal
bundle $\IP(N_{V/X})$.

\[
  \begin{tikzcd}
    E \arrow[r, "j", hook] \arrow[d, "p"'] & \widetilde X \arrow[d, "\pi"] \\
    V \arrow[r, "i"', hook] & X
  \end{tikzcd}\label{intersection-squar-2}\tag{$\ast$}
\]

Let $T$ be an irreducible subvariety of $X$ of codimension $d$, and assume
that $T\not\subset V$. We denote by $\widetilde T$ the proper transform of $T$,
that is, the Zariski closure of $\pi^{-1}(T\setminus V)$ in $\widetilde X$.
By the blow-up closure lemma \cite{ravi2022rising} or
\cite[\href{https://stacks.math.columbia.edu/tag/080E}{Tag 080E}]{stacks-project}, the proper
transform $\widetilde T$ is canonically identified with the blow-up of $T$ along the subscheme
$T \cap V$ and has pure dimension $m-d$. Thus we have the following commutative diagram.

\[
  \begin{tikzcd}[row sep=1em, column sep =1.5em]
    \widetilde T\cap E \arrow{dd}\arrow{rr} \arrow{dr} &  & \widetilde T \arrow{dd}\arrow{dr} &  &  & \\
    {} & E \arrow[dd, crossing over] \arrow[rr, crossing over] &  & \widetilde X \arrow{dd} &  & \\
    T\cap V \arrow{dr}\arrow{rr} &  & T \arrow{dr} &  &  & \\
    {} & V\arrow[rr, crossing over] &  & X &  & \\
  \end{tikzcd}\label{intersection-cubic}\tag{$\diamond$}
\]

Since $E$ is a hypersurface in a smooth variety which does not contain $\widetilde T$, the
intersection $\widetilde T\cap E$ is proper. It follows that $\widetilde T\cap E$ has pure dimension
$m-d-1$, hence defines an effective cycle of codimension $d$ inside $E$. The fiber of
$\widetilde T\cap E\to T\cap V$ over a point $p\in T\cap V$ may have any dimension between $0$ and
$e-1$. When $T$ and $V$ intersect properly, every fiber of $\widetilde T\cap E\to T\cap V$ has
dimension $e-1$, and is in fact a projective space $\IP^{e-1}$.

By Theorem~\ref{structure-of-chow-ring-of-projective-bundle},
the fundamental class of $\widetilde T\cap E$ can be written uniquely in the form
\begin{equation}\label{chow-ring-structure-projective-bundle}
  [\widetilde T \cap E]
  =
  \sum_{\ell=0}^{e-1}p^\ast(a_\ell)h^\ell,
  \qquad
  a_\ell\in A^{d-\ell}(V).
\end{equation}
We use the convention $A^q(V)=0$ if $q<0$ or $q>\dim V$.

Let $C = C_{T\cap V}T$ be the normal cone of $T\cap V$ in $T$ and $\IP(C)$ its projectivization.
Then there is a canonical isomorphism
$$
  \IP(C) \simeq \widetilde T\cap E.
$$
Let $\iota_T\colon T\cap V\hookrightarrow V$ be the inclusion. Then we have
\[
  s_V(T\cap V,T)
  \coloneqq
  (\iota_T)_*s(T\cap V,T) = p_\ast\left( \sum_{i\geq 0} h^i\cdot[\widetilde T\cap E ] \right)\in A_\ast(V) .
\]
In order to measure the norm of a Chow class, we may identify the class in $A^\ast(V)$ with its
image in the numerical space $N^*_\R(V)$.

\begin{proposition}\label{proposition-the-three-equivalence-volume-estimates}
  Fix an integer $d$, a norm on the spaces $N^\ast_\R(V)$, and an ample line
  bundle $\widetilde L$ on $\widetilde X$. The following uniform
  estimates are equivalent:

  \begin{enumerate}
    \item There is a constant $C >0$ such that for all $T\sst X$ irreducible subvarieties of
      codimension $d$ not contained in $V$,
      \[
        \max_{0\leq\ell\leq e-1}\|a_\ell\|
        \leq C\Vol_L(T).
      \]
    \item There is a constant $C>0$ such that for all $T\sst X$ irreducible subvarieties of
      codimension $d$ not contained in $V$,
      \[
        \|s_V(T\cap V,T)\|
        \leq C \Vol_L(T).
      \]
    \item There is a constant $C >0$ such that for all $T\sst X$ irreducible subvarieties of
      codimension $d$ not contained in $V$,
      \[
        \Vol_{\widetilde L}(\widetilde T)
        \leq C\Vol_L(T).
      \]
  \end{enumerate}
  The three constants need not be equal and may depend on the fixed data, but
  not on $T$.
\end{proposition}

\begin{proof}
  Put $R_T=[\widetilde T\cap E]$. The projective bundle formula identifies
  $R_T$ with its coefficient vector $(a_0,\ldots,a_{e-1})$ via
  \eqref{chow-ring-structure-projective-bundle}
  \[
    R_T=\sum_{\ell=0}^{e-1}p^\ast(a_\ell)h^\ell\in N^{d}_\IR(E) .
  \]
  One calculates
  \[
    p_*(h^q R_T)
    =\sum_{\ell+q\geq e-1}
    a_\ell\,s_{\ell+q-e+1}\in N^{d+q-e+1}_\IR(V) ,
    \qquad q\geq0,
  \]
  where $s_{\ell+q-e+1} = s_{\ell+q-e+1}(N_{V/X})\coloneqq p_\ast(h^{\ell+q})$ denotes the
  corresponding Segre class of $N_{V/X}$.

  Write
  $S_j\coloneqq \{s_V(T\cap V,T)\}_{m-e-j}\in N^{j}_\IR(V)$.
  By definition, $S_j=0$ for $j<d-e+1$, whereas, for $j\geq d-e+1$,
  \[
    S_j
    =p_*(h^{j-d+e-1}R_T)
    =\sum_{\ell=\max\{0,d-j\}}^{e-1}
    a_\ell\,s_{\ell+j-d}.
  \]

  Explicitly, the transformation from $(a_\ell)$ to $(S_{d-r})$ is described by the upper-triangular
  matrix
  \[
    \begin{pmatrix}
      S_d\\
      S_{d-1}\\
      \vdots\\
      S_{d-e+2}\\
      S_{d-e+1}
    \end{pmatrix}
    =
    \begin{pmatrix}
      1 & s_1 & s_2 & \cdots & s_{e-1}\\
      0 & 1 & s_1 & \cdots & s_{e-2}\\
      \vdots & \ddots & \ddots & \ddots & \vdots\\
      0 & \cdots & 0 & 1 & s_1\\
      0 & \cdots & \cdots & 0 & 1
    \end{pmatrix}
    \begin{pmatrix}
      a_0\\
      a_1\\
      \vdots\\
      a_{e-2}\\
      a_{e-1}
    \end{pmatrix},
    \qquad s_k=s_k(N_{V/X}),
  \]
  This matrix is upper triangular with identity diagonal. By induction, one can prove
  $$
    \|s_V(T\cap V,T)\|
    \leq e\sum_{k=0}^{e-1}||s_k|| \max_{0\leq\ell\leq e-1}\|a_\ell\|
  $$
  and
  \[
    \max_{0\leq\ell\leq e-1}\|a_\ell\|\leq \sum_{k=0}^{e-1}||s_k||\, \|s_V(T\cap V,T)\|.
  \]
  Therefore (1) and (2) are equivalent.

  By \eqref{norm-and-volume-comparison}, for sufficiently small $\epsilon >0$, the line bundle
  $\widetilde L = \pi^\ast L -\epsilon E$ is ample
  \cite[Proposition~1.7.10]{lazarsfeld2017positivityI}. Expanding
  $c_1(\widetilde L)^{m-d}\cdot[\widetilde T]$, we have
  \[
    \pi^\ast\left( c_1(L)^{m-d}\cdot[T] \right)
    +j_\ast\left( \sum_{q=1}^{m-d}\binom{m-d}{q}(-\epsilon)^q
    (j^\ast\pi^\ast L)^{m-d-q}(j^*E)^{q-1}\cdot R_T\right).
  \]
  The second term is
  \[
    -j_\ast p^\ast \left( \sum_{q=1}^{m-d}\binom{m-d}{q} \epsilon^q
    (i^*L)^{m-d-q}\, S_{d+q-e} \right).
  \]
  It is a fixed linear combination of $S_{d-e+1},S_{d-e+2},\ldots, S_d$.
  The first term is exactly $\Vol_{L}(T)$, so \eqref{norm-and-volume-comparison} shows that (2) and
  (3) are equivalent.
\end{proof}

\begin{theorem}\label{complexity-inequality}
  Fix a norm on the spaces $N^\ast_\R(V)$. There exists a constant
  $C>0$, depending only on $(X,V,L)$. For every irreducible subvariety $T\subset X$ not contained in
  $V$ and $a_\ell \in A^{d-\ell}(V)$ satisfying
  \[
    [\widetilde T \cap E]
    =
    \sum_{\ell=0}^{e-1}p^\ast(a_\ell)h^{\ell},
  \]
  we have
  \[
    \|a_\ell\|\leq C\Vol_L(T).
  \]
\end{theorem}

For convenience, we reverse the indexing of the sum,
\[
  [\widetilde T \cap E]
  =
  \sum_{\ell=0}^{e-1}p^\ast(b_\ell)h^{e-1-\ell},
  \qquad b_\ell=a_{e-1-\ell}\in A^{d-e+1+\ell}(V).\label{b_j} \tag{$\circledast$}
\]
We first establish several useful lemmas.

\begin{lemma}\label{lemma-complexity-lower}
  For each $\ell\in\{0,\ldots,e-1\}$ there are coefficients $p_j\in A^{l-1-j}(V)$ independent of
  $T$ such that
  \[
    b_\ell+ \sum_{j=0}^{\ell-1} p_j\cdot b_j
  \]
  is pseudo-effective in $A^{d-e+1+\ell}(V)$, where $b_j$ are defined by \eqref{b_j}.
\end{lemma}

\begin{proof}
  Set $\omega=c_1(L|_V)$. By Lemma~\ref{lemma-relative-ample}, choose an integer
  $C_0>0$ so that $h+C_0p^\ast\omega$ is ample on $E$. The class
  \[
    p_\ast\bigl((h+C_0p^\ast\omega)^{\ell}
    [\widetilde T\cap E]\bigr)
  \]
  is pseudo-effective. The projective-bundle formula gives
  \[
    p_\ast\bigl((h+C_0p^\ast\omega)^{\ell}
    [\widetilde T\cap E]\bigr)
    =b_{\ell}+B_{\ell},\qquad
    B_{\ell}=\sum_{j=0}^{\ell-1} b_j\,
    p_\ast\bigl(h^{e-1-j}(h+C_0p^\ast\omega)^{\ell}\bigr).
  \]
  Terms with $j>\ell$ vanish under $p_\ast$, and the coefficient of
  $b_{\ell}$ is $p_\ast(h^{e-1})=1$. Each coefficient
  \[
    p_j=p_\ast\bigl(h^{e-1-j}(h+C_0p^\ast\omega)^{\ell}\bigr)
  \]
  is independent on $T$.
\end{proof}

\begin{lemma}\label{lemma-complexity-hypersurfaces}
  There is an integer $M>0$, such that for every
  irreducible subvariety $T\subset X$ not contained in $V$ one can find hyperplane sections
  $H_1,\ldots,H_e\in|ML|$ containing $V$ with $\pi^\ast H_i=\widetilde H_i+E$ such that for every
  $1\leq r\leq e$, both $H_1\cap\cdots\cap H_r$ and
  $\widetilde T\cap\widetilde H_1\cap\cdots\cap\widetilde H_r$
  are proper intersections.
\end{lemma}

\begin{proof}
  Choose $M\in \IN$ large enough, independently of $T$, such that $\mathcal I_V \otimes M L$ is
  globally generated and
  $A_M = M \pi^\ast L + \mathcal O_{\widetilde X}(-E)$
  is very ample. By
  $\pi_\ast\mathcal O_{\widetilde X}(-E)=\mathcal I_V$, the space
  $$
    \Lambda = H^0(X,\mathcal I_V \otimes M L)\simeq H^0(\widetilde X,A_M )
  $$
  gives a base-point-free system on
  $U = X\setminus V$ and the complete system of $A_M$ on
  $\widetilde X$. So every section $s\in \Lambda$ gives a hypersurface $H_s$ containing $V$ on $X$
  and a hypersurface $D_s$ on $\widetilde X$.

  Then we choose $s_1,\ldots,s_e\in \Lambda$ successively. For each $1\leq r\leq e-1$, suppose that
  $H_1,\ldots, H_r$ are chosen such that $H_1\cap\cdots\cap H_r$ and
  $\widetilde T\cap\widetilde H_1\cap\cdots\cap\widetilde H_r$
  are proper. We apply Bertini's theorem in arbitrary characteristic
  \cite[Theorem~3.2]{jelonek2013bertini} to every irreducible component $C$ of
  $H_1\cap\cdots\cap H_r\cap U$. For a general section $s \in \Lambda$, the zero set of $s|_C$ is
  empty or has dimension $\dim C-1$. Then we check such a section properly cuts the irreducible
  component of $H_1\cap\cdots\cap H_r$ containing $V$. Let $C$ be one of these irreducible
  components. Then the zero locus of $s|_{C\setminus V}$ is empty or has dimension $\dim C -1$ and
  $s$ vanishes along $V$. Since $\dim C > \dim V$, we have $H_{s}$ cuts $C$ properly.

  We apply the same Bertini theorem to every irreducible component of
  $\widetilde T\cap\widetilde H_1\cap\cdots\cap\widetilde H_r$. For a general section
  $s\in \Lambda$, the intersection $D_s\cap C$ is either empty or has dimension $\dim C-1$. Thus
  $\widetilde T\cap\widetilde H_1\cap\cdots\cap\widetilde H_r\cap D_s$ is proper.

  Since $A_M$ is very ample, we may also require the divisor $D_s$ on $\widetilde X$ not to contain
  $E$. Then $D_s=\widetilde H_s$, and the defining identity
  $\mathcal I_V\mathcal O_{\widetilde X}=\mathcal O_{\widetilde X}(-E)$
  gives $\pi^\ast H_s=\widetilde H_s+E$.

  Finally, let $H_{r+1}=H_s$. We find the $r+1$-th desired hypersurface. We can conclude the lemma
  by induction.
\end{proof}

\begin{lemma}\label{lemma-complexity-upper}
  Fix $\ell\in\{0,\ldots,e-1\}$. There exists hyperplanes
  $H_1,\ldots,H_r\in|ML|$ as in Lemma~\ref{lemma-complexity-hypersurfaces}.
  There are coefficients $q_{j}\in N^{\ell-j}_\IR(V)$ independent on $T$ and \(H_i\), such that
  \[
    i_\ast b_{\ell}
    =T\cdot H_1\cdots H_{\ell+1}
    -\pi_\ast(\widetilde T\cdot\widetilde H_1\cdots\widetilde H_{\ell+1})
    +i_\ast\left(\sum_{j=0}^{\ell-1} q_j b_j \right).
  \]
  where $i_\ast$ is the push-forward in \eqref{intersection-squar-2}.
\end{lemma}

\begin{proof}
  Let $\omega=c_1(L|_V)$ and $r=\ell+1$. we recall the commutative diagram
  \eqref{intersection-squar-2}:
  \[
    \begin{tikzcd}
      E \arrow[r, "j", hook] \arrow[d, "p"'] & \widetilde X \arrow[d, "\pi"] \\
      V \arrow[r, "i"', hook] & X
    \end{tikzcd}
  \]
  On $E$ we have
  $j^\ast E=-h$ and $j^\ast\pi^\ast H_i = p^\ast i^\ast H_i=M p^\ast\omega$. Hence
  $j^\ast\widetilde H_i=h+M p^\ast \omega$. By \eqref{chow-ring-structure-projective-bundle}, we
  have $j^\ast( \widetilde T )= [\widetilde T\cap E]= \sum b_\ell h^{e-1-\ell}$.

  Now consider the intersection
  \begin{align*}
    \Xi\coloneqq T\cdot H_1 \cdots H_r = \pi_\ast\left( \widetilde T\cdot \pi^\ast H_1\cdots \pi^\ast H_r \right).
  \end{align*}
  Expanding $\pi^\ast H_1=\widetilde H_1+E$ one factor, we get
  \begin{align*}
    \Xi = \pi_\ast\left( \widetilde T\cdot \widetilde H_1\cdot\pi^\ast H_2\cdots \pi^\ast H_r \right)
    + \pi_\ast\left( \widetilde T\cdot E\cdot\pi^\ast H_2\cdots \pi^\ast H_r \right).
  \end{align*}
  Continuing to expand one factor at a time, we obtain
  \begin{align*}
    \Xi = \pi_\ast\left(\widetilde T\cdot \widetilde H_1\cdots \widetilde H_r\right) +\sum_{t=1}^{r} \pi_\ast\left( \widetilde T\cdot \widetilde H_1\cdots \widetilde H_{t-1} \cdot E\cdot \pi^\ast H_{t+1}\cdots\pi^\ast H_r\right)
  \end{align*}
  Applying the projection formula and expanding in $A^\ast(E)$, the $t$-th term in the previous
  formula equals
  \begin{align*}
    \pi_\ast & \left( j_\ast [\widetilde T\cap E] \cdot \widetilde H_1 \cdots \widetilde H_{t-1} \cdot \pi^\ast H_{t+1}\cdots\pi^\ast H_r \right)\\
    & = \pi_\ast j_\ast\left( [\widetilde T\cap E]\cdot j^\ast \left( \widetilde H_1\cdots \widetilde H_{t-1}\cdot\pi^\ast H_{t+1}\cdots\pi^\ast H_r\right)\right) \\
    & =i_\ast p_\ast \left(\sum_{\ell=0}^{e-1} p^\ast(b_\ell) h^{e-1-\ell}\cdot \left(h+ M p^\ast\omega \right)^{t-1} \cdot (Mp^\ast\omega)^{r-t} \right).
  \end{align*}
  Thus we obtain
  \begin{align*}
    \Xi = \pi_\ast & \left(\widetilde T\cdot \widetilde H_1\cdots \widetilde H_r\right)\\
    & + i_\ast p_\ast\left( \left(\sum_{\ell=0}^{e-1} p^\ast(b_\ell) h^{e-1-\ell}\right) P_r(h,M\omega) \right),
  \end{align*}
  where $P_r$ is a polynomial in $h$ and $x$,
  \[
    P_r(h,x)=\sum_{t=1}^r(h+x)^{t-1}x^{r-t}
    =\sum_{t=0}^{r-1}\binom r{t+1}h^t x^{r-1-t}.
  \]

  By Lemma~\ref{lemmaprojectionformulaforprojectivebundle}, the class
  $p_\ast\left( h^{e-1-j}P_r(h,x) \right)$ vanishes when $j>r-1$ and we have
  $p_\ast\left( h^{e-1-j}P_r(h,x) \right) = 1$ when $j=r-1$. For $j<r-1$, set
  \[
    c_{\ell,j}
    =p_\ast\bigl(h^{e-1-j}P_r(h,Mp^\ast\omega)\bigr)
    \in N^{\ell-j}_\R(V).
  \]
  These are fixed classes, expressed in terms of $\omega$ and the Segre
  classes of $N_{V/X}$. Then we obtain what we want
  \begin{align*}
    T\cdot H_1 \cdots H_r = \pi_\ast & \left(\widetilde T\cdot \widetilde H_1\cdots \widetilde H_r\right) + i_\ast\left(b_{\ell} + \sum_{j=0}^{\ell-1}b_j \cdot c_{\ell,j} \right).
  \end{align*}
\end{proof}

\begin{proof}[Proof of Theorem~\ref{complexity-inequality}]
  We argue by induction on $\ell$ that $||b_\ell||\leq C\Vol_L(T)$ for some constant $C$. We firstly
  consider the base case $b_0$. If
  $d-e+1+\ell$ lies outside $[0,\dim V]$, then $b_{0}=0$. Otherwise,
  Lemma~\ref{lemma-complexity-upper} shows $b_0$ is pseudo-effective. We apply
  Lemma~\ref{lemma-complexity-upper} to $\ell =0$ and obtain
  \[
    i_\ast b_0 \leq T\cdot H_1.
  \]
  Thus $||b_0|| \leq M\Vol_L(T)$.

  For $\ell>1$, assume that the desired estimate is known
  for every $b_j$ with $j<\ell$. For convenience, put $r=\ell+1$ and
  $k=\dim V-q=m-d-r\geq0$.

  Lemma~\ref{lemma-complexity-lower} supplies a pseudo-effective class
  $$
    \beta=b_{\ell}+B_{\ell}
  $$
  for $B_\ell=\sum_{j=0}^{\ell-1} p_j\cdot b_j$. The induction hypothesis gives
  \[
    \|B_{\ell}\|\leq C_1\Vol_L(T).
  \]
  By Lemma~\ref{lemma-complexity-upper}, there exist hypersurfaces $H_1,\ldots,H_{\ell+1}$ as in
  Lemma~\ref{lemma-complexity-hypersurfaces} and we obtain \(A_{\ell}=\sum_{j=0}^{\ell-1}q_j b_j\)
  such that
  \[
    i_\ast b_{\ell}=T\cdot H_1\cdots H_r-\pi_\ast(\widetilde T\cdot\widetilde H_1\cdots\widetilde H_r)+i_\ast A_{\ell}.
  \]
  We have $\|A_{\ell}\| \leq C_2 \Vol_L(T)$ where $C_2$ depends only on $(X,V,L)$.

  The proper intersection implies $\pi_\ast(\widetilde T\cdot\widetilde H_1\cdots\widetilde H_r)$ is
  effective. Intersecting this
  identity with $c_1(L)^k$ and using the projection formula yields
  \begin{align*}
    b_{\ell}\cdot\omega^k
    & =(T\cdot H_1\cdots H_r-\pi_\ast(\widetilde T\cdot\widetilde H_1\cdots\widetilde H_r)+i_\ast A_{\ell})\cdot c_1(L)^k\\
    & \leq M^r\Vol_L(T)+i_\ast( A_{\ell}\cdot c_1(L|_V)^k).
  \end{align*}
  It leads to \(b_{\ell}\cdot\omega^k\leq C_3 \Vol_L(T)\) for an independent constant \(C_3\).

  The norm--volume comparison~\eqref{norm-and-volume-comparison} on $V$ gives a constant $C'$. Then,
  we have
  \[
    \|b_{\ell}\|
    \leq\|\beta\|+\|B_{\ell}\|\leq C'\,\left( \beta\cdot \omega^k \right) +\|B_{\ell}\|.
  \]
  Combining these inequalities, we conclude
  \begin{align*}
    \|b_{\ell}\|
    \leq & C'C_3\Vol_L(T) + C_1 \Vol_L(T).
  \end{align*}

  Taking an independent constant \(C = C'C_3+C_1\). This proves the induction step. Since there are
  only finitely many
  degrees $d$ and indices $\ell$, the constants can be chosen
  uniformly.
\end{proof}

\subsection{Corollaries of volume inequalities}

\begin{theorem}\label{theorem-complexity-inequality-under-birational-morphism}
  Assume that $\varphi\colon Y\to X$ is a birational morphism between smooth projective varieties
  over a field of characteristic $0$. Fix ample line bundles $L_X\to X$ and $L_Y\to Y$.

  Then there exists a constant $C>0$ such that, if $T\subset X$ is a
  subvariety not contained in $\Ind(\varphi^{-1})$ and
  $\widetilde T$ denotes its proper transform under $\varphi$, then
  \begin{align*}
    \Vol_{ L_Y}(\widetilde T)\leq C \Vol_{L_X}(T).
  \end{align*}
\end{theorem}

We refer to \cite[Lemma~2.5]{vu2022exotic} for an analytic version.

\begin{proof}
  Suppose first that $\varphi$ is a composition of blow-ups along smooth
  centers. The assertion then follows immediately from
  Theorem~\ref{complexity-inequality} and
  Proposition~\ref{proposition-the-three-equivalence-volume-estimates}.

  We now consider an arbitrary birational morphism
  $\varphi\colon Y\to X$. Since $\varphi$ is projective and birational, it is
  the blow-up of a coherent ideal sheaf $\cI\subseteq\cO_X$; see, for
  instance, \cite[Theorem~II.7.17]{hartshorne2013algebraic}.

  By \cite[Theorem~19]{kollar2005resolution}, there exists a composition of
  blow-ups along smooth centers
  \[
    f\colon X'\longrightarrow X
  \]
  such that $X'$ is smooth and $\cI\cO_{X'}$ is invertible. The universal
  property of the blow-up therefore gives a birational morphism
  $\psi\colon X'\to Y$ satisfying $f=\varphi\circ\psi$.

  Let $T'\subset X'$ be the proper
  transform of $T$ under $f$. Then
  $\psi(T')=\widetilde T$. Choose an ample line bundle $L_{X'}$ on $X'$.
  For some integer $M\gg 0$, the line bundle
  \[
    H\coloneqq ML_{X'}-\psi^*L_Y
  \]
  is ample.

  Put $r=\dim T=\dim T'$. Since $T' \subset X'$ is not
  contracted by $\psi$, the projection formula and the nefness of
  $\psi^*L_Y$ give
  \[
    \Vol_{L_Y}(\widetilde T)
    \leq \Vol_{\psi^*L_Y}(T')
    \leq M^r\Vol_{L_{X'}}(T').
  \]
  Since $f$ is a composition of blow-ups at smooth centers, we conclude that
  \[
    \Vol_{L_{X'}}(T') \leq C\Vol_{L_X}(T).
  \]
  for some constant $C$.
\end{proof}

\bigskip

Let $X$ be a smooth projective variety and fix an ample line bundle $L$ on
$X$. Given a subscheme $Z$ of $X$, we set $n(Z)\in\IN$ to be the number of irreducible components of
$Z$.

\begin{theorem}\label{theorem-number-of-intersection-irred-components}
  There exists a constant $C>0$ such that for any two irreducible subvarieties
  $Z_1,Z_2\subset X$,
  \[
    n(Z_1\cap Z_2) \leq C\Vol_L(Z_1)\Vol_L(Z_2),
  \]
  where $Z_1\cap Z_2$ is the scheme-theoretic intersection.

\end{theorem}

We first prove the following lemma:

\begin{lemma}
  Under the assumptions and notation of Theorem~\ref{complexity-inequality}, there exists a constant
  $C>0$, independent of $T$, such that
  \[
    n(T\cap V) \leq C\Vol_L(T).
  \]

\end{lemma}
\begin{proof}
  Fix an ample line bundle $A$ on $E$. We consider \(D_1,\cdots, D_r\) all the irreducible component
  of closed subscheme $\widetilde T\cap E$. Set \(D_i' = \pi(D_i)\), then \(D_i'\) are irreducible
  closed subset covering \(T\cap V = \pi(\widetilde T\cap E)\). Therefore,
  \[
    n(T\cap V) \leq r=n(\widetilde T\cap E) .
  \]

  The $A$-volume of each irreducible component is at least $1$. Hence
  \[
    n(T\cap V)
    \leq \Vol_A(\widetilde T\cap E)
    \leq C\Vol_L(T),
  \]
  where the last inequality follows from Theorem~\ref{complexity-inequality} and
  Proposition~\ref{proposition-the-three-equivalence-volume-estimates}.
\end{proof}

\begin{proof}[Proof of Theorem~\ref{theorem-number-of-intersection-irred-components}]
  When $Z_1$ and $Z_2$ are of dimension $0$, the theorem is obvious. For the remaining case, we have
  $Z_1\times Z_2 \not\subset \Delta$.
  Apply the preceding lemma to $X\times X$, to the diagonal
  $\Delta\subset X\times X$, and to $T=Z_1\times Z_2$. With
  $\overline L=\pi_1^\ast L+\pi_2^\ast L$, one has
  \[
    \Vol_{\overline L}(Z_1\times Z_2)
    =
    \binom{\dim Z_1+\dim Z_2}{\dim Z_1}
    \Vol_L(Z_1)\Vol_L(Z_2).
  \]
  The binomial coefficient is bounded in terms of $X$, and the result follows.
\end{proof}

\begin{remark}
  If we replace $e(Z)$ by the number of irreducible components with geometric multiplicities, the
  theorem does not hold anymore.
\end{remark}

\section{Periodic points and multiplicity}\label{section-periodic-pts}
Let $K= \overline K$ be an algebraically closed field.
\subsection{Isolated periodic points and multiplicity}
Let $X$ be a normal projective variety of dimension $m$ defined over $K$, and let
$f\colon X\dasharrow X$ be a dominant rational self-map. Recall that $f$ is determined by its graph
$\Gamma(f)\subset X\times X$ so that the first projection $\pi_1\colon \Gamma(f) \to X$ is
a birational morphism. The locus where the inverse $\pi_1^{-1}$ is not well-defined is the
indeterminacy locus $\Ind(f) \subset X$: it is a closed subset of codimension at least $2$. The
graph-fixed scheme of $f^n$ is $F_n(f) \coloneqq \Gamma(f^n)\times_{X\times X} \Delta$

\begin{definition}
  The scheme of periodic locus of period \(n\) is
  \[
    \Per_n(f) = \pi_1(F_n(f)).
  \]
\end{definition}

In general, the indeterminacy locus is not empty and the set of periodic points has infinite
cardinality, so it is natural to consider the set of isolated periodic points.

Let $U_n$ be the open locus on which the length-$n$ orbit under $f$ is well defined. More precisely,
we can define $U_n$ by induction,
$$
  U_0 =X,\qquad U_{n+1} = \left\{x\in U_n, f^n(x)\not\in \Ind(f) \right\}.
$$

\begin{definition}

  The scheme of isolated (genuine) periodic point is
  \[
    \IPer_n(f) = \Per_n(f)^{\text{iso}} \times_X U_n,
  \]
  where \(\Per_n(f)^{\text{iso}}\) is the subscheme of \(0\)-dimensional irreducible components.

  Throughout, ``period \(n\)''
  means fixed by \(f^n\), rather than necessarily having least period
  \(n\).

\end{definition}

\begin{definition}
  We introduce two counts in \(\IPer_n(f)\).

  \begin{enumerate}
    \item The number of isolated periodic points without multiplicity, denoted by
      \[
        \NIPer_n(f) \coloneqq \#\left(\IPer_n(f)(K) \right) = \deg_K(\IPer_n(f)_{\text{red}}).
      \]

    \item The number of isolated periodic points with multiplicity is denoted by
      \[
        \NIPerM_n(f)\coloneqq
        \sum_{x\in\operatorname{IPer}_n(f)} l_x(f^n) =\deg_K(\IPer_n(f))
      \]
      where \(l_x(f^n)\) denotes the geometric multiplicity of $P= (x,x)$ in $F_n(f)$; see
      \eqref{geometric-multiplicity} for the definition.
  \end{enumerate}

\end{definition}

\begin{remark}
  \begin{enumerate}
    \item One could count with Samuel multiplicities which is
      \[
        e_x(f^n)\coloneqq \left(e_{F_n(f)}(\Gamma(f^n) )\right)_P
      \]
      the multiplicity of $\Gamma(f^n)$ along scheme $F_n(f)$ at irreducible component $P=(x,x)$;
      see Section~\ref{section-Multiplicity along subvariety} for details.
    \item When \(X\) is smooth, the multiplicity $e_x(f^n)$ is, in fact, the intersection
      multiplicity $i(P, \Gamma(f^n)\cdot\Delta;X\times X)$. Thus we may use the notation
      $i_x(f^n) = e_x(f^n)$ for intersection multiplicity. Then the counts with geometric multiplicity
      and intersection multiplicity are in fact same, because $f^n$ is regular near $x$, so
      $A=\cO_{\Gamma(f^n),P}$ is regular and hence Cohen--Macaulay; see
      \cite[Proposition~11.1.10]{swanson2006integral}.
  \end{enumerate}
\end{remark}

\subsection{Lefschetz number and perfect points}

we assume further that $X$ is smooth. Note that the diagonal \(\Delta \hookrightarrow X\times X\) is
a regular embedding if and only if \(X\) is smooth.

Then one can use intersection theory to define the \textit{Lefschetz number} of $f^n$, that is the
intersection number
\[
  L(f^n)\coloneqq\Gamma(f^n)\cdot \Delta \in \Z~.
\]

The growth of Lefschetz numbers could relate to cohomological data (étale cohomological data in
positive characteristic) via Lefschetz fixed-point theorem.

In general, we do not have easy comparison between $L(f^n)$ and $\NIPer_n(f)$ as well as other
counts with or without multiplicity. Because the higher-dimensional component and multiplicities may
contribute wildly to the intersection number.

Following Fulton's approach, we consider the blow-up of the diagonal
$\pi\colon \widetilde{X\times X} \to X\times X$ and look at the schematic intersection of
the exceptional divisor $E\coloneqq\pi^{-1}(\Delta)$ with the strict transform
$\widetilde{\Gamma(f^n)}$ of the graph.
The crucial observation is that
$E$ is a hypersurface, so the intersection $E\cap \widetilde{\Gamma(f^n)}$
has pure dimension $m-1$.

\begin{definition}
  We call $P =(x,x) \in \Delta(K)$ a perfect point of period $n$ if the fiber
  $E\cap \widetilde{\Gamma(f^n)} \to \Delta$ over $P$ has full dimension i.e. contains $E_P$ the
  full fiber of $E$ over $P$.

  Denote by $\Perfect_n(f)$ the number of perfect points of $\Gamma(f^n)$ with geometric
  multiplicity of $E_P$ in $E\cap \widetilde{\Gamma(f^n)}$.
\end{definition}

\begin{lemma}\label{lemma-NGG-P}
  All isolated points are perfect and we have
  \[
    \NIPer_n(f)\leq \NIPerM_n(f) \leq \Perfect_n(f).
  \]
\end{lemma}
\begin{proof}
  Put $\Gamma=\Gamma(f^n)$, $Z=\Gamma\cap\Delta$ and $R=E\cap\widetilde{\Gamma}$, and let $P=(x,x)$
  be an isolated point of $Z$. The strict transform $\widetilde\Gamma$ is the blow-up of $\Gamma$
  along $Z$, and its exceptional divisor is
  \[
    R=\IP(C_Z\Gamma);
  \]
  see \cite[Section~16.2]{fulton2013intersection}. The normal cone $C_Z\Gamma$ is pure of dimension
  $m$ and embeds in $N_{\Delta/(X\times X)}|_Z$. Since $P$ is isolated, its support over $P$ is
  therefore the whole vector space $N_{\Delta/(X\times X)}|_P\simeq\mathbb A^m$. Projectivizing
  shows that $E_P\simeq\IP^{m-1}$ is an irreducible component of $R$, so $P$ is perfect.

  By \cite[Definition~7.1]{fulton2013intersection}, $i_x(f^n)$ is the coefficient of this component
  in $[C_Z\Gamma]$. Projectivization preserves its generic multiplicity, hence
  \[
    \ell(E_P,R)=i_x(f^n).
  \]
  Thus each isolated graph-periodic point contributes its intersection multiplicity to
  $\Perfect_n(f)$. All other perfect points have positive multiplicity, so summing gives
  $\NIPerM_n(f)\leq \Perfect_n(f)$.
\end{proof}

\begin{table}[ht]
  \centering
  \small
  \renewcommand{\arraystretch}{1.3}
  \begin{tabular}{p{0.57\textwidth}p{0.20\textwidth}p{0.12\textwidth}}
    \hline
    Description & Notation & Type \\
    \hline
    Periodic locus
      & \(\Per_n(f)\) & Scheme \\
    Isolated genuine periodic points
      & \(\IPer_n(f)\) & Scheme \\
    Isolated periodic points, counted without multiplicity
      & \(\NIPer_n(f)\) & Count \\
    Isolated periodic points, counted with geometric multiplicity
      & \(\NIPerM_n(f)\) & Count \\
    Lefschetz number
      & \(L(f^n)\) & Count \\
    Perfect points, counted with geometric multiplicity
      & \(\Perfect_n(f)\) & Count \\
    \hline
  \end{tabular}
  \caption{Periodic-point notation.}
\end{table}

\subsection{Under birational conjugacy}
Let $X$ be a normal complex projective surface and $p\in X$. Let $\pi_p\colon Y \to X$ be the
blow-up at $p$. For any dominant rational map $f\colon X\dasharrow X$, set
$g =\pi_p^{-1}\circ f\circ \pi_p\colon Y\dasharrow Y$ the lift of $f$. We show that periodic-point
counts of \(f\) can be bounded in terms of the counts of \(g\). This result will be used in
Section~\ref{section3}.

\begin{lemma}\label{lemma-birationa-conjugacy-upper-bound}
  There exists a constant $C =C(f,\pi_p)$ such that, for every $n\in \IN$, we have
  \[
    \NIPerM_n(f) - \NIPerM_n(g) \leq C(f,\pi_p).
  \]
\end{lemma}

\begin{proof}
  Every isolated periodic point of $f$ whose orbit avoids $p$ lift to an isolated periodic point of
  $g$ with the same multiplicity.

  If $p$ is an isolated periodic point of $f$, let $r$ be its exact period. Then
  $\NIPerM_n(f) \leq \NIPerM_n(g)$ when $r\nmid n$. When $n = rm$ for some $m\in \N$, we have
  \[
    \NIPerM_n(f) - \NIPerM_n(g) \leq r\, l_p(f^n).
  \]
  since the orbit of $p$ has $r$ points.

  By Lemma~\ref{boundedness-of-multiplicity}, we have \(l_p(f^n)\) is bounded. This concludes the
  proof.
\end{proof}

When \(p\) is a smooth point on \(X\). The boundedness of multiplicity is shown by
\cite{shub1974remark}. We give a proof in the case when \(p\) is a singular point following the
local calculation in \cite{shub1974remark}.
\begin{lemma}\label{boundedness-of-multiplicity}
  Let \(h\colon (Z,p)\to(Z,p)\) be a holomorphic self-map of a complex analytic germ and \(p\) is an
  isolated fixed point. The multiplicity \(l_p(h^n)\) is bounded.
\end{lemma}
\begin{proof}
  Embed \((Z,p)\) into \((\IC^N,0)\) and extend \(h\) to a holomorphic germ
  \(H\colon (\IC^N,0)\to (\IC^N,0)\). Let \(A = DH(0)\) and set \(M\) as the least common multiple
  of the orders of the root-of-unity eigenvalues of \(A\).

  For \(n\geq 1\), set \(k=\gcd(n,M)\). Then \(n=sk\) for \(s\in \IN\). We claim that
  \begin{align}
    H^n(z)-z = Q(z)\left(H^k(z)-z\right)\label{H^n}
  \end{align}
  for \(Q(z)\) an invertible holomorphic matrix germ. Indeed, one can write
  \[
    H^n(z)-z= \sum_{t=0}^{s-1} H^{tk}\left(H^{k}(z)\right)- H^{tk} \left( z\right)
  \]
  with \(H^{tk}\left(w\right)- H^{tk} \left( z\right) =B_t(w,z)(w-z)\) and
  \(B_{t}(0,0) =DH^{tk}(0).\)
  Thus
  \[
    Q(z) = \sum_{t=0}^{s-1} B_t(H^k(z) ,z), \qquad Q(0) = \sum_{t=0}^{s-1}DH^{tk}(0) = \sum_{t=0}^{s-1} A^{tk}.
  \]
  Let \(\alpha\) be any eigenvalue of \(A\), then
  \(q(\alpha)\coloneqq 1+\alpha^{k}+\cdots+\alpha^{(s-1)k}\) is an eigenvalue of \(Q(0)\). When
  \(\alpha^k =1\), then \(q(\alpha) =s \not =0\). When \(\alpha^k\not=1\), we have
  \(\alpha^{n}=\alpha^{sk}\not=1\). Indeed, let \(r\) be the order of \(\alpha\), if \(r|n\) then
  \(r|\gcd(n,M)=k\) which is a contradiction. Then
  \(q(\alpha) =\frac{1-\alpha^n}{1-\alpha^k}\not=0\). We conclude that \(Q(0)\) is invertible. Thus
  \(Q\) is invertible.

  Let \(J_k\subset\cO_{Z,a}\) be the fixed ideal of \(h^k\) i.e. the local ring at \(p\) of scheme
  \(\pi_1\left(F_k(h) \right)\). By \eqref{H^n}, we have ideals
  \[
    \left(H^n(z)-z \right) = \left(H^k(z)-z \right)\sst \cO_{\IC^N,0}.
  \]
  Restricting to \(Z\) gives \(\cJ_n=\cJ_k\). We can conclude that
  \[
    l_p(h^n) = l_p\left(h^{\gcd(n,M)}\right),
  \]
  and they are bounded by \(\max\{l_p(h^{k}), k|M\}\).

\end{proof}

\section{Upper bound on the number of periodic points}\label{section3}

\subsection{Rational maps} \label{section-upper-bound-rational-map}
Let $X$ be a smooth projective variety of dimension $m$ over $K=\overline K$, and let $L$ be an
ample line bundle on $X$. Let $\pi_1$ and $\pi_2$ be the projections onto the first and second
factors.

We consider the blow-up of $X\times X$ along the diagonal $\Delta$:
\[
  \begin{tikzcd}
    E \arrow[r, "j", hook] \arrow[d, "p"'] & \widetilde{X\times X} \arrow[d, "\pi"] \\
    \Delta \arrow[r, "i"', hook] & X\times X
  \end{tikzcd}
\]
The exceptional divisor $E$ is a projective bundle over $\Delta$ and is isomorphic to the
projectivized tangent bundle $\IP T_X$.
Let $f\colon X\dasharrow X$ be a dominant rational self-map of $X$. Let $\widetilde{\Gamma(f)}$
denote the strict transform of the graph ${\Gamma(f)}$ under the blow-up.

\begin{proof}[Proof of Theorem~\ref{theoremA}]
  We fix an ample line bundle $\widetilde L$ on $\widetilde{X\times X}$ and consider the ample line
  bundle $\overline{L} = \pi_1^\ast L + \pi_2^\ast L$ on $X\times X$.

  Let $P(f)$ be the number of perfect points of $\Gamma(f)\cap \Delta$ counted with geometric
  multiplicity. Let $P=(x,x)\in \Delta(K)$ be a perfect point. Therefore, we have the following
  inequality of effective cycles:
  $$
    [\widetilde{\Gamma(f)} \cap E]\geq P(f) [E_{P}].
  $$
  Therefore, we obtain
  \begin{align*}
    \Vol_{\widetilde L }\left(P(f) [E_{P}]\right) \leq \Vol_{\widetilde L}([\widetilde{\Gamma(f)} \cap E]).
  \end{align*}

  Theorem~\ref{complexity-inequality} and
  Proposition~\ref{proposition-the-three-equivalence-volume-estimates}, applied to
  $(X \times X, \Delta, \overline{L}, \widetilde L)$, yield a constant $C_1$ such that
  \begin{align*}
    P(f) & \leq \Vol_{\widetilde L}([\widetilde{\Gamma(f)} \cap E]) \leq C_1 \cdot \Vol_{\overline{L}}({\Gamma(f)}).
  \end{align*}
  The volume of the graph ${\Gamma(f)}$ with respect to $\overline{L}$ is given by
  \begin{align*}
    \Vol_{\overline{L}}({\Gamma(f)}) = {\Gamma(f)} \cdot \left(\pi_1^\ast L + \pi_2^\ast L \right)^m = \sum_{p=0}^m \binom{m}{p} \deg_{L,p}(f)\leq 2^m \max_{0 \leq p \leq m} \deg_{L,p}(f).
  \end{align*}
  It follows that
  $$
    P(f) \leq 2^m C_1 \max_{0 \leq p \leq m} \deg_{L,p}(f).
  $$
  Replacing $f$ with its $n$-th iterate $f^n$, we deduce
  \begin{align*}
    \Perfect_n(f)
    \leq C \max_{0 \leq p \leq m} \lambda_p(f^n).
  \end{align*}
  where $C =2^m C_1$. This completes the proof.
\end{proof}

\subsection{Birational surface maps}\label{section-birational-surface-map}

We now prove the upper-bound part of Theorem~\ref{theoremB}. In fact, we will prove the following
theorem which is stronger and gives an upper bound on $\NIPerM_n(f)$.

\begin{theorem}\label{upper-bound-surface}
  Let $X$ be a normal projective surface over an algebraically closed field of characteristic $0$.
  Then, for any birational map $f\colon X\dasharrow X$ with $\lambda_1(f) >1$, we are in one of the
  following two cases:
  \begin{enumerate}
    \item The surface $X$ is birationally isomorphic to an abelian surface and
      \[
        \NIPerM_n(f) \leq \lambda_1(f)^n + O\bigl(\lambda_1(f)^{n/2}\bigr).
      \]
    \item The surface $X$ is not birationally isomorphic to an abelian surface and
      \[
        \NIPerM_n(f) \leq \lambda_1(f)^n + O(1).
      \]
  \end{enumerate}
\end{theorem}
Set $\lambda_1=\lambda_1(f)$. Choose a finitely generated subfield
$k_0/\IQ$ over which $X$ and $f$ are defined. By embedding its algebraic closure into $\IC$, we may
assume that $X$ and $f$ are defined over $\IC$.

We claim that to prove Theorem~\ref{upper-bound-surface}, it suffices to prove it under the
assumption that $X$ is smooth and $f$ is algebraically stable.

To achieve this assumption, we take a resolution of singularities of $X$ and then apply
\cite[Theorem~0.1]{diller2001dynamics}. We obtain
a birational morphism $\pi\colon\widetilde X\to X$ for a smooth projective variety $\widetilde X$
and an algebraically
stable lift $\widetilde f\colon\widetilde X\dasharrow\widetilde X$. By
Lemma~\ref{lemma-birationa-conjugacy-upper-bound}, there is a constant $C=C(f,\pi)$ such that
\[
  \NIPerM_n(f)\leq \NIPerM_n(\widetilde f)+C.
\]
Also the first dynamical degree is invariant under birational conjugacy. So it suffices to prove the
theorem under hypothesis of smoothness and algebraic stability.

We shall compare the Lefschetz number $\Gamma(f^n)\cdot \Delta$ with the number $\NIPer_n(f)$ of
isolated periodic points. We first estimate the Lefschetz numbers.

\begin{lemma}\label{lemma-growth-of-lefschetz-number}
  Let \(X\) be a smooth complex projective surface and $f: X\dasharrow X$ an algebraically stable
  birational map. We have
  \begin{enumerate}
    \item When the surface $X$ is birationally equivalent to an abelian surface, then
      \[
        \Gamma(f^n)\cdot \Delta = \lambda_1^n + O\bigl(\lambda_1^{n/2}\bigr).
      \]
    \item If the surface $X$ is not birationally equivalent to an abelian surface, then
      \[
        \Gamma(f^n)\cdot \Delta = \lambda_1^n + O(1).
      \]
  \end{enumerate}
\end{lemma}

\begin{proof}
  By the Lefschetz fixed-point formula, for every $n\geq 1$ we have
  \begin{align*}
    \Gamma(f^n) \cdot \Delta = \sum_{i=0}^{4} (-1)^i\operatorname{Trace}\left((f^n)^\ast\colon H^i(X,\IC) \to H^i(X,\IC) \right).
  \end{align*}
  Algebraic stability gives $(f^n)^\ast=(f^\ast)^n$ on $H^{1,1}(X)$.
  By \cite[Theorem~5.1]{diller2001dynamics}, the spectrum of $f^\ast$ on
  $H^{1,1}(X)$ outside the unit disk consists of the single simple eigenvalue
  $\lambda_1$. On the other hand, resolving the graph shows that birational
  pull-back preserves the Hodge norm on $H^{2,0}(X)$ and $H^{0,2}(X)$, so all
  eigenvalues on these two spaces have modulus $1$. It follows that
  \[
    \operatorname{Trace}\left((f^n)^\ast\colon H^2(X,\IC) \to H^2(X,\IC) \right) = \lambda_1^n + O(1).
  \]

  By \cite{cantat2001autom}\cite[Theorem~B]{blanc2016dynamical}, $X$ is birationally equivalent to a
  rational surface, an abelian surface, a K3 surface, or an Enriques surface.
  In the first, third, and fourth cases, $H^1(X,\IC)=H^3(X,\IC)=0$. In the
  abelian case, these two cohomology groups are birationally invariant and may
  be computed on the abelian minimal model. The induced automorphism there has
  spectral radius $\lambda_1^{1/2}$ on both $H^1$ and $H^3$. Therefore
  \[
    \left|\operatorname{Trace}\left((f^n)^\ast\colon
    H^1(X,\IC)\oplus H^3(X,\IC)\right)\right|
    =O\bigl(\lambda_1^{n/2}\bigr).
  \]
  Finally, the traces on $H^0(X,\IC)$ and $H^4(X,\IC)$ are both $1$.
  Combining this with the preceding estimates gives
  \[
    \Gamma(f^n)\cdot\Delta=
    \begin{cases}
      \lambda_1^n+O\bigl(\lambda_1^{n/2}\bigr),
      & \text{if $X$ is birational to an abelian surface},\\
      \lambda_1^n+O(1), & \text{otherwise}.
    \end{cases}
  \]
\end{proof}

Then we investigate the following example which allows us to calculate the vertical and horizontal
contributions to the intersection number.
\begin{example}\label{example-4-dim}
  Consider a $4$-dimensional smooth variety $X$ together with a $2$-dimensional smooth subvariety
  $V\subset X$. We consider the blow-up along $V$:
  \[
    \begin{tikzcd}
      E \arrow[r, "j", hook] \arrow[d, "p"'] & \widetilde X \arrow[d, "\pi"] \\
      V \arrow[r, "i"', hook] & X
    \end{tikzcd}
  \]
  where $N =\cN_{V,X}$ is the rank-$2$ normal bundle of $V$ and $E = \IP N$ is a $\IP^1$-bundle.

  Let $T$ be a possibly singular subvariety of dimension 2 of $X$ not containing $V$, and let
  $\widetilde T$ be its strict transform.

  Let $C$ be an irreducible component of $\widetilde T\cap E$ with geometric multiplicity
  $\ell_C\geq 1$. We describe the contribution of $C$ to the intersection number $(T\cdot V)$,
  denoted by $(T\cdot V)^C$. By Definition~6.1.2 and Example~6.1.8 of \cite{fulton2013intersection},
  we can write
  \begin{align*}
    T\cdot V = \sum (T\cdot V)^C,
  \end{align*}
  where $(T\cdot V)^C = \ell_C\, \{p_\ast(c_{1}(\xi)\cdot [C])\}_{0}$ and $\xi= p^\ast N/ \cO_E(-1)$
  is the universal quotient line bundle.

  We say $C$ is a vertical component, if $Z = p(C)$ has dimension $0$. Then $C =p^{-1}(Z)$ is a
  fiber of the projective bundle. Thus, we have
  \[
    (T\cdot V)^C =\ell_C\in \IN.
  \]

  We say $C$ is a horizontal component, if $Z=p(C)$ has dimension $1$. Write $p_*[C]=k_C[Z]$ with
  $k_C\in\IN^\ast$. Since $c_1(\xi)=p^*c_1(N)+h$, where
  $h=c_1(\cO_E(1))$, the projection formula gives
  \[
    (T\cdot V)^C
    =\ell_C\bigl(k_C[Z]\cdot c_1(N)+[C]\cdot h\bigr).
  \]
\end{example}

\begin{proposition}\label{proposition-horizontal-lift}
  Let $X$ be a smooth complex projective surface, let
  $f\colon X\dasharrow X$ be birational with $\lambda_1(f)>1$, and let
  $Z\subset X$ be an irreducible curve that is fixed pointwise by $f$ away from
  $\Ind(f)$. Blow up the diagonal
  of $X\times X$, and denote the exceptional divisor by $E$ and its projection
  to the diagonal by $p$. For every $n\geq 1$, there is a unique irreducible
  curve
  \[
    C_n\subset \widetilde{\Gamma(f^n)}\cap E
  \]
  such that $p(C_n)=Z$. Moreover, the geometric multiplicities
  $\ell_{C_n}(f^n)\coloneqq
  \ell(C_n,\widetilde{\Gamma(f^n)}\cap E)$ are bounded independently of $n$, and
  $p_*[C_n]=[Z]$.
\end{proposition}

\begin{proof}[Proof of Theorem~\ref{upper-bound-surface}]

  Let $\Delta\colon X\hookrightarrow X\times X$ be the diagonal embedding. We
  blow up along the diagonal and denote by $E$ the exceptional divisor, by
  $p\colon E\to\Delta\simeq X$ the natural projection, and by
  $\widetilde{\Gamma(f^n)}$ the strict transform of $\Gamma(f^n)$.

  We now analyze the intersection along $E$. Let $\cC_n$ be the finite set of
  irreducible components of $\widetilde{\Gamma(f^n)}\cap E$. For a horizontal component $C\in\cC_n$,
  set $Z = p(C)\subset\Delta\simeq X$. Then $Z$ is an irreducible curve fixed pointwise by $f^n$
  away
  from $\Ind(f^n)$. By \cite{cantat2010invariant} and \cite[Theorem~5.5]{xie2012periodic}, there are
  only finitely many
  periodic curves; hence the set
  \[
    \cZ = \{ Z = p(C) \mid C\in \bigcup_{n\geq 1} \cC_n \text{ and $C$ is horizontal} \}
  \]
  is finite.

  Put $N=N_{\Delta/(X\times X)}\simeq T_X$ and
  $h=c_1(\cO_E(1))$. According to Example~\ref{example-4-dim}, the
  contribution of $C$ to $\Gamma(f^n)\cdot\Delta$ is
  \[
    (\Gamma(f^n)\cdot \Delta)^C = \ell_C(f^n) \left( k_C\, [Z]\cdot c_1(N) + m_C\right)
  \]
  where $\ell_C(f^n)=\ell(C,\widetilde{\Gamma(f^n)}\cap E)$ is the
  geometric multiplicity, $p_*[C]=k_C[Z]$ with $k_C\in\IN^\ast$, and
  \[
    m_C=[C]\cdot h\in\Z.
  \]

  One can write for each $n\in \N$
  \begin{align*}
    \Gamma(f^n)\cdot \Delta = \sum_{C\in \cC_n, \dim p(C)=0}(\Gamma(f^n)\cdot \Delta)^C + \sum_{Z\in \cZ}\sum_{C\in \cC_n, p(C) =Z}(\Gamma(f^n)\cdot \Delta)^C = A_n + B_n.
  \end{align*}

  Here $A_n$ is the contribution of the vertical components and $B_n$ is the
  contribution of the horizontal components. The former is exactly the count of perfect points with
  multiplicity. In particular, by Lemma~\ref{lemma-NGG-P}, we have
  \[
    \NIPerM_n(f)\leq \Perfect_n = A_n.
  \]
  It remains to show that $B_n$ is uniformly bounded from below. For each
  $Z\in\cZ$, let $q_Z$ be its least pointwise period. If $q_Z\nmid n$, no
  horizontal component of $\widetilde{\Gamma(f^n)}\cap E$ lies over $Z$. If
  $n=q_Zm$, Proposition~\ref{proposition-horizontal-lift}, applied to
  $f^{q_Z}$, gives a unique component $C_{Z,n}$ over $Z$ such that
  $p_*[C_{Z,n}]=[Z]$, and its multiplicity is bounded independently of $n$.
  Thus
  \[
    B_n=\sum_{\substack{Z\in\cZ\\q_Z\mid n}}
    \ell_{C_{Z,n}}(f^n)
    \bigl([Z]\cdot c_1(N)+[C_{Z,n}]\cdot h\bigr).
  \]
  Choose an ample line bundle $L$ on $\Delta$, put $\omega=c_1(L)$, and
  choose $A\gg0$ so that $h+A p^*\omega$ is ample on $E$. Then
  \[
    [C_{Z,n}]\cdot h>-A[C_{Z,n}]\cdot p^*\omega
    =-A[Z]\cdot\omega.
  \]
  Since $\cZ$ is finite and the multiplicities above are uniformly bounded,
  $B_n$ is uniformly bounded from below. Consequently,
  \[
    \NIPerM_n(f)\leq \Gamma(f^n)\cdot\Delta-B_n\leq\Gamma(f^n)\cdot\Delta+O(1).
  \]

  Combining this estimate with Lemma~\ref{lemma-growth-of-lefschetz-number} proves the theorem.
\end{proof}

\begin{proof}[Proof of Proposition~\ref{proposition-horizontal-lift}]
  Consider a general smooth point $x\in Z\setminus\Ind(f)$, where $f$ is a
  local biholomorphism, and choose analytic coordinates $(z_1,z_2)$ centered at \(x\) with
  $Z=\{z_1=0\}$. Since $f$ fixes $Z$ pointwise, there is a largest integer
  $t\geq1$ such that
  \[
    f(z_1,z_2)=
    \bigl(z_1+z_1^tg(z_1,z_2),z_2+z_1^th(z_1,z_2)\bigr),
  \]
  where at least one of $g,h$ is a unit at a general point of $Z$. We define
  $\nu_Z(f)=t$.

  Near \((x,x)\in X\times X\) the blow-up of the diagonal is given by
  \[
    \bigl\{(z,w,[l_0:l_1])\in\IC^4\times\IP^1
    \mid (z_1-w_1)l_1=(z_2-w_2)l_0\bigr\}.
  \]
  Suppose first that $g$ is a unit, and use the chart \(\{l_0 \not =0\}\) with
  \[
    (z_1,z_2,a,u) \mapsto (z_1,z_2,z_1+a, z_2+au, [1:u]).
  \]
  The strict transform of the graph is given by the ideal
  \[
    (a -z_1^t g, u - h/g).
  \]
  Its intersection with $E =\{a=0 \}$ is given by the ideal
  \[
    (z_1^t,u-h/g).
  \]
  This is a subscheme with unique irreducible component with multiplicity \(t\). We conclude that
  there is a unique irreducible component \(C\) and its multiplicity is $\nu_Z(f)$. We have further
  that \(C\) is given by
  \[
    z_1=0,\qquad a=0,\qquad u=\frac{h(0,z_2)}{g(0,z_2)}.
  \]
  And the projection \(p\) deduces an isomorphism
  \[
    \left(0,z_2,0,\frac{h(0,z_2)}{g(0,z_2)}\right) \longmapsto (0,z_2).
  \]
  So \(p_\ast[C] =[Z]\).

  The chart $\{l_1 \not=0\}$ gives the same conclusion when $h$ is a unit.

  Applying this calculation to every iterate yields irreducible curve \(C_n\) with multiplicity
  \(\ell_{C_n}(f^n) = \nu_Z(f^n)\) and \(p_\ast[C_n] =[Z].\)

  It remains to bound $\nu_Z(f^n)$. For any smooth point $z=(0,z_2) \in Z\backslash\Ind(f)$, the
  differential of \(f\) has the form
  \[
    df|_{z}=\begin{pmatrix}\lambda(z_2) &0\\ \mu(z_2) &1\end{pmatrix},
    \qquad \lambda(z_2)\in \IC(Z)\backslash \{0\} ,\quad \mu(z_2)\in \IC(Z).
  \]

  One may write
  \[
    f(z_1,z_2) = \left(\lambda(z_2) z_1 + O(z_1^2), z_2 + \mu(z_2)z_1 + O(z_1^2) \right).
  \]

  If $\lambda(z_2)$ is not a root of unity for some \(z_2\), then $\nu_Z(f^n) =1$ for all
  \(n\in \IN\). Indeed, the $n$-th iterate $f^n$ has the form
  \begin{align*}
    f^n(z_1,z_2) =\left( \lambda(z_2)^nz_1+ O(z_1^2), z_2+ \mu(z_2)\left(\sum_{j=0}^{n-1} \lambda(z_2)^j \right)z_1 +O(z_1^2) \right).
  \end{align*}

  Thus we can suppose that \(\lambda\) is a constant root of unity, then $\nu_Z(f^n) = 1$ whenever
  $\lambda^n\neq 1$. Hence, by replacing $f$ by an iterate, we may assume $\lambda = 1$. If
  $\mu\not=0$, we have $\nu_Z(f^n) = 1$.

  So we may assume further $\mu=0$, then $\nu_Z(f^n) \geq 2$. Then we rewrite \(f\) in the following
  form
  \[
    f(z_1,z_2) =(z_1+P(z_1,z_2),z_2+Q(z_1,z_2)) = (z_1+\sum_{n\geq 2 } P_n(z_1,z_2), z_2+\sum_{n\geq 2 } Q_n(z_1,z_2)),
  \]
  where \(P,Q\) are formal polynomials and $P_n, Q_n$ are homogeneous polynomials of degree $n$. By
  assumption, $\nu_Z(f)=\min\{ \ord_{z_1}(P), \ord_{z_1}(Q)\} \geq 2$.

  By \cite{brochero2008parabolic}, there exists a unique formal vector field
  \[
    \cX = a(z_1,z_2)\frac{\partial}{\partial z_1} + b(z_1,z_2)\frac{\partial}{\partial z_2} = \sum_{n\geq 2}a_n(z_1,z_2)\frac{\partial}{\partial z_1} + \sum_{n\geq 2}b_n(z_1,z_2)\frac{\partial}{\partial z_2},
  \]
  such that
  \[
    f(z_1,z_2) = (\exp(\cX)(z_1),\exp(\cX)(z_2)),
  \]
  where $a_n, b_n$ are homogeneous polynomials of degree $n$ and
  $\exp(\cX)(g) = \sum_{j=0}^{+\infty} \frac{1}{j!}\cX^j(g)$ for $g\in \IC[[z_1,z_2]]$.
  For $q\geq 2$, we set
  $\cX_q = \sum_{2\leq n\leq q}a_n(z_1,z_2)\frac{\partial}{\partial z_1} + \sum_{2\leq n\leq q}b_n(z_1,z_2)\frac{\partial}{\partial z_2}$.
  Then we have \(P_2 = a_2\) as well as \(Q_2 =b_2\) and for \(q\geq 2\),
  \begin{align}\label{parabolic-eq1}
    P_{q+1} = a_{q+1} + \operatorname{HT}_{q+1}\left( \sum_{j=2}^{q}\frac{1}{j!}\cX_q^j(z_1) \right),\\
    Q_{q+1} = b_{q+1} + \operatorname{HT}_{q+1}\left( \sum_{j=2}^{q}\frac{1}{j!}\cX_q^j(z_2) \right),\label{parabolic-eq2}
  \end{align}
  where $\operatorname{HT}_{q+1}(h)$ denotes the homogeneous term of $h$ of order $q + 1$.
  \smallskip

  We claim that
  \[
    \min\{\ord_{z_1}(P),\ord_{z_1}(Q)\}
    =
    \min\{\ord_{z_1}(a),\ord_{z_1}(b)\}.
  \]
  Set \(t=\min\{\ord_{z_1}(P),\ord_{z_1}(Q)\}\). Then \(z_1^t\) divides both \(P_q\) and \(Q_q\) for
  every \(q\). In particular, \(z_1^t\) divides \(a_2\) and \(b_2\). Using \eqref{parabolic-eq1} and
  \eqref{parabolic-eq2}, we obtain inductively that \(z_1^t\) divides \(a_n\) and \(b_n\) for every
  \(n\). Conversely, set \(s=\min\{\ord_{z_1}(a),\ord_{z_1}(b)\}\). Then \(z_1^s\) divides \(a_n\)
  and \(b_n\) for every \(n\). Equations \eqref{parabolic-eq1} and \eqref{parabolic-eq2} imply that
  \(z_1^s\) divides \(P_2\) and \(Q_2\), and then, inductively, \(P_{q+1}\) and \(Q_{q+1}\) for
  every \(q\).

  Hence,
  \[
    \nu_Z(f) =\min \{\ord_{z_1}(P),\ord_{z_1}(Q) \} = \min \{\ord_{z_1}(a),\ord_{z_1}(b)\}.
  \]
  Then
  \[
    \nu_Z(f^n) =\nu_Z(n\cX) = \min \{\ord_{z_1}(na),\ord_{z_1}(nb)\}= \nu_Z(f).
  \]
  We conclude that \(C_n\) has bounded multiplicity \(\ell_{C^n}(f^n)=\nu_Z(f^n)\).
\end{proof}

\section{Moriwaki heights and the lower bound}
\label{section4}

In this section, let $K$ be a finitely generated field over $\IQ$. Let $k$ be the algebraic closure
of $\Q$ inside $K$. Then $k$ is a finite extension of $\IQ$. For details, see
\cite{moriwaki2000arithmetic}, \cite{chen2020arakelov}, and \cite{zhang1995small}.

\subsection{Product formula}\label{section4-1}

Let $K$ be a finitely generated field over $\IQ$ of transcendence degree $b$. Let \(\mathscr B\) be
a normal, projective, flat integral scheme over \(\operatorname{Spec}\mathbf Z\), with function
field \(K\) and smooth generic fiber. Fix an arithmetically ample \(C^\infty\)-Hermitian line bundle
\(\overline{\mathscr H}\) on \(\mathscr B\), whose curvature form is positive on
\(\mathscr B(\mathbf C)\). We fix such a model $(\sB,\overline{\sH})$ of $K$.

Let $\Omega_K^\infty$ denote the set of all embeddings of $K$ into $\IC$. One may view
$\Omega_K^\infty$ as a subset of $\sB(\IC)$ via the composition $\Spec(\IC)\to\Spec(K)\to\sB$. Then
\begin{align*}
  \Omega_K^\infty = \sB(\IC)\setminus \bigcup_{Y\in \Div(\sB)}Y(\IC),
\end{align*}
where $Y$ runs over all divisors on $\sB$ defined over $\IQ$. For each $\sigma\in\Omega_K^\infty$,
the associated place is defined by
\[
  |\cdot|_\sigma \coloneqq |\sigma(\cdot)|_\IC.
\]

Let $\Omega_K^{\mathrm{NA}}$ be the set of all prime divisors on $\sB$. Any
$\Gamma\in\Omega_K^{\mathrm{NA}}$ defines a non-archimedean absolute value via
\[
  |\cdot|_\Gamma \coloneqq \exp\bigl(-\ord_\Gamma(\cdot)\bigr).
\]

Note that horizontal divisors give non-archimedean places trivial on $k$, while vertical divisors
extend non-archimedean places of $k$.

The set of places of $K$, denoted by $\Omega_K$ (or simply $\Omega$), is the disjoint union
\[
  \Omega_K \coloneqq \Omega_K^\infty\sqcup \Omega_K^{\mathrm{NA}}.
\]
We equip $\Omega_K^{\mathrm{NA}}$ with the discrete topology and $\Omega_K^\infty$ with the analytic
topology.

On $\Omega_K$, we consider the measure
\begin{align*}
  \nu = \sum_{\Gamma}\delta_{\Gamma}\cdot \widehat{\deg}(\overline{\sH}^b|_\Gamma) + \operatorname{MA}(c_1(\overline{\sH})),
\end{align*}
where $\operatorname{MA}(c_1(\overline{\sH})) = c_1(\overline{\sH})^{\wedge b}$ is the
Monge--Amp\`ere measure on $\sB(\IC)$. Since $\overline{\sH}$ is ample, the measure $\nu$ is
nonnegative. The product formula then holds: for any $a\in K\setminus\{0\}$,
\begin{align*}
  \int_{\Omega_K} \log|a|_w\,d\nu(w)
  = \sum_{\Gamma} \log|a|_\Gamma\, \widehat{\deg}(\overline{\sH}^b|_\Gamma)
  + \int_{\sB(\IC)} \log|a|_w\,\operatorname{MA}(c_1(\overline{\sH})) = 0.
\end{align*}

We say that a function $\varphi\colon \Omega_K\to\IR$ is $\nu$-\textit{dominated} if there exists a
$\nu$-integrable function $\widetilde\varphi$ such that $|\varphi|\le \widetilde\varphi$.

\subsection{Strongly adelic metrics}

Let $X$ be a projective variety over $K$, and let $w\in\Omega_K$ be a place of $K$. Let $K_w$ be the
corresponding completion and $X_{K_w}$ the base change. Let $X_w$ denote the associated analytic
space: $X_w = X(\overline{K_w})$ when $w$ is archimedean, and $X_w$ is the Berkovich analytic space
of $X_{K_w}$ when $w$ is non-archimedean. For a line bundle $L$ on $X$, denote by $L_w$ its
analytification on $X_w$.

\begin{definition}
  Let $L$ be a line bundle on $X$. A \textit{metric family} on $L$ is a family
  $\varphi=\{\|\cdot\|_w\}_{w\in\Omega_K}$ of continuous metrics on the induced line bundles over
  the analytic spaces $X_w$.
\end{definition}

\begin{definition}
  A $\sB$-model of $(X,L)$ consists of the data $(\sX, \sL,\|\cdot\| )$: $\sX$ is a projective
  arithmetic variety over $\sB$, $\sL$ is a line bundle on $\sX$ such that $\sX_K=X$ and $\sL_K=L$,
  and $\|\cdot\|$ is a continuous semipositive metric on $\sL_\IC$ over $\sX(\IC)$.

  A metric family is called model if it is induced by a model of $(X,L)$.
\end{definition}

\begin{definition}
  We say that a metric family $\varphi= \{\|\cdot\|_w\}_{w\in \Omega_K}$ on $L$ is
  \textit{FS (Fubini--Study)} if, at every place $w\in \Omega_K$, the metric $\|\cdot\|_w$ has the
  form
  \begin{align*}
    & \|s\|_{w}(x) = \frac{|s(x)|}{\max( |\tau_1(x)|,|\tau_2(x)|,\ldots, |\tau_s(x)| )^{1/q}},
  \end{align*}
  where $q\geq 1$ is an integer and the $\tau_i$ are finitely many sections of $qL$ without common
  zeros.
\end{definition}

\begin{definition}
  We say that a metric family $\varphi$ on $L$ is \textit{strongly adelic} if there exist
  $(M,\varphi_1)$ and $(N,\varphi_2)$ such that $L=M-N$, $\varphi=\varphi_1-\varphi_2$, and both
  $\varphi_1$ and $\varphi_2$ are model and FS metric families.
\end{definition}

\begin{remark}
  Every strongly adelic metric family is measurable and dominated in the sense of Chen–Moriwaki
  \cite[Section~6.1]{chen2020arakelov}.
\end{remark}

\begin{example}
  Let $L$ be a base-point-free line bundle on $X$. Let $\Sigma = \{\tau_1,\ldots, \tau_s \}$ be a
  finite subset generating $H^0(X,L)$. Then $\Sigma$ defines a morphism from $X$ to some projective
  space. Hence $\Sigma$ defines an FS metric family on $L$.

  In fact, one can prove that this metric family is also a model metric family. Hence such a metric
  family is strongly adelic.
\end{example}

\begin{definition}
  If $\varphi=\{\|\cdot\|_w\}_{w\in\Omega_K}$ and $\varphi'=\{\|\cdot\|_{w}'\}_{w\in\Omega_K}$ are
  two strongly adelic metric families on $L$, then $\varphi-\varphi'$ defines a metric family on
  $L\otimes L^{\vee}=\cO_X$. We define the local distance by
  \begin{align*}
    d_w(\varphi,\varphi') = \sup_{x\in X_w}\left|\log\|\cdot\|_w(x)-\log\|\cdot\|_{w}'(x)\right|,
  \end{align*}
  and the global distance between $\varphi$ and $\varphi'$ is defined as the integral
  \begin{align*}
    \dist(\varphi,\varphi') = \int_{w\in\Omega_K} d_w(\varphi,\varphi')\,d\nu(w) \in [0,+\infty].
  \end{align*}

\end{definition}

\subsection{Moriwaki height}

Let $\overline{L}$ be a line bundle equipped with an adelic FS model metric. For $x\in X(K)$, the
Moriwaki height is
\begin{align*}
  h_{\overline{L}}(x) = \int_{\Omega_K}-\log\|s(x)\|_w d\nu(w),
\end{align*}
where $s$ is a section of $L$ not vanishing at $x$.

The integral is well defined because the metric is induced by a model and the height can be
interpreted in terms of arithmetic intersection numbers.

Denote by $h_{\overline{L},s,w}(x)=-\log\|s(x)\|_w$ the local height of $\overline{L}$ at the place
$w$.

\begin{example}
  The coordinate sections define the naive height function for $(\IP^N, \cO(1))$. We set
  \begin{align*}
    h_{\mathrm{nv}}([x_0:\cdots:x_N]) \coloneqq \int_{\Omega_K} \log\max(|x_0|_w,\ldots,|x_N|_w)\,d\nu(w),
  \end{align*}
  for $[x_0:\cdots:x_N]\in \IP^N(K)$.
\end{example}

\begin{definition}[Distance on $X_w$]
  We fix an embedding $\iota:X\hookrightarrow \IP^N_K$. For every place $w\in\Omega_K$, we define a
  distance $\dist_w$ on the analytic space $X_w$ by pulling back the spherical distance on
  $\IP^N(K_w)$ (which is associated with the standard FS metric). More explicitly, the spherical
  distance is defined by the formula
  \begin{align*}
    \dist_w([x_0:\cdots:x_N], [y_0:\cdots:y_N])=\frac{\max_{0\leq i<j\leq N} |x_iy_j-x_jy_i|_w}{\max_{0\leq i\leq N}|x_i|_w\,\max_{0\leq j\leq N}|y_j|_w}.
  \end{align*}

\end{definition}

\begin{proposition}\label{proposition-boundedness-effective-divisor-height}
  Let $X$ be a projective variety defined over $K$. Let $D\in \Div(X)$ be an effective divisor
  equipped with a strongly adelic metric $\varphi$ on $\cO(D)$. Then the height $h_{(D,\varphi )}$
  is bounded from below on $(X\backslash \supp(D))(K)$.
\end{proposition}

\begin{proof}

  First, we can construct a strongly adelic metric $\varphi_0$ as follows. Write $D=D_1-D_2$, where
  $D_1$ and $D_2$ are very ample line bundles. Choose a basis $\Sigma=\{f_1,\ldots,f_s\}$ of
  $\Gamma(X,\cO(D_2))$. Then for any $i$ we have
  \begin{align*}
    \div(f_i) +D_1 = D+\div(f_i)+D_2\geq D \geq0.
  \end{align*}
  Therefore $f_i\in\Gamma(X,\cO(D_1))$. Hence we may extend $\Sigma$ to a basis
  $\Sigma'=\{f_1,\ldots,f_s,\ldots,f_t\}$ of $\Gamma(X,\cO(D_1))$. The model metric families induced
  by $\Sigma$ and $\Sigma'$ on $\cO(D_2)$ and $\cO(D_1)$ define an adelic metric $\varphi_0$ on
  $\cO(D)$ by tensor product.

  With respect to $\varphi_0$, we check that the local heights are nonnegative: for any
  $w\in \Omega_K$ and $x\not\in \supp(D_1)(K)$,
  \begin{align*}
    h_{(D,\varphi_0) ,w}(x) & = h_{\overline{D}_1,w}(x)-h_{\overline{D}_2,w}(x)\\
    &= \log\max\{ |f_1(x)|_w,\ldots, |f_t(x)|_w \} -\log\max\{ |f_1(x)|_w,\ldots, |f_s(x)|_w \}\\
    &\geq 0.
  \end{align*}
  Since the measure $\nu$ is nonnegative, the Moriwaki height
  \begin{align*}
    h_{(D,\varphi_0)}(x) = \int_{\Omega_K} h_{(D,\varphi_0) ,w}(x) d\nu(w)\geq 0.
  \end{align*}

  For any other strongly adelic metric $\varphi$ on $\cO(D)$, Proposition 6.1.5 of
  \cite{chen2020arakelov} implies that the local distance function
  $(w\in\Omega_K)\mapsto d_w(\varphi,\varphi_0)$ is $\nu$-dominated. Thus
  \[
    h_{(D,\varphi)}(x) \ge h_{(D,\varphi_0)}(x)-\int_{\Omega_K} d_w(\varphi,\varphi_0)\,d\nu(w),
  \]
  so $h_{(D,\varphi)}(x)$ is bounded from below for $x\not\in \supp(D_1)(K)$.

  In order to obtain a uniform lower bound on $(X\backslash\supp(D))(K)$, we consider a very ample
  line bundle $A$ such that $A+\cO(D)$ is also very ample. Let $g_1,\cdots, g_r$ be a basis of
  $H^0(X,A)$. Then each $g_i$ gives a decomposition $D = \left( \div(g_i)+D \right) -\div(g_i)$. We
  repeat the discussion above, and have a uniform lower bound on
  \[
    \bigcup_{i=1}^r
    \bigl(X\setminus\supp(\div(g_i)+D )\bigr)(K)= \bigl(X\setminus\supp(D)\bigr)(K).
  \]

\end{proof}

\subsection{Differences of heights}
Let $X$ be a smooth projective surface defined over $K$ and let $L$ be a line bundle on $X$. Let
$f\colon X\dasharrow X$ be a dominant rational map defined over $K$. Let $f = g\circ \pi^{-1}$ be
the minimal resolution of $f$. Note that $f^\ast L = \pi_\ast g^\ast L$.
\[
  \begin{tikzcd}
    \widetilde{X} \arrow[d, "\pi"'] \arrow[rrd, "g"] &  & \\
    X \arrow[rr, "f", dashed] &  & X
  \end{tikzcd}
\]

\begin{proposition}\label{prop-difference-height-nef-divisor}
  Suppose that $D$ is a nef divisor. For any strongly adelic metrics on $\cO(D)$ and
  $\cO(f^\ast D)$, we have
  \begin{align*}
    h_{\overline{D}}(f(x))\leq h_{\overline{f^\ast D}}(x)+O(1),
  \end{align*}
  for $x\not \in \Ind(f)$ and $x\not\in \supp(D)$.
\end{proposition}

\begin{proof}
  We take the minimal smooth resolution where $\pi\colon \widetilde X\to X$ is a composition of
  blow-ups and $g\colon \widetilde X \to X$ is the lift of $f$.

  By \cite{diller2001dynamics}, $E=\pi^\ast f^\ast D-g^\ast D$ is effective. Moreover,
  $\overline{E}=\pi^\ast(\overline{f^\ast D})-g^\ast\overline{D}$ carries an adelic metric. Let $y$
  be the unique preimage of $x$ under $\pi$.
  \begin{align*}
    h_{\overline{D}}(f(x))- h_{\overline{f^\ast D}}(x)= h_{\overline{D}}(g(y))-h_{\overline{f^\ast D}}(\pi(y)) = h_{g^\ast \overline{D}}(y)-h_{\pi^\ast \overline{f^\ast D}}(y) =- h_{\overline{E}}(y).
  \end{align*}
  Since $E$ is effective, this yields the desired upper bound.

\end{proof}

Let $U=X\setminus\Ind(f)$. Since $X$ is smooth and $\Ind(f)$ has codimension at least $2$, there is
a natural identification
\[
  H^0(X,f^\ast L)\xrightarrow{\sim}H^0(U,f^\ast L|_U).
\]
Then for any section $\sigma\in H^0(X,L)$, we define $T_f(\sigma)\in H^0(X,f^\ast L)$ to be the
unique extension of $(f|_U)^\ast \sigma$. This defines a morphism
\[
  T_f: H^0(X,L) \to H^0(X,f^\ast L),
\]
and $T_f(\sigma)/T_f(\tau) = (\sigma/\tau) \circ f$ as rational functions. By the construction, we
have also
\[
  \div(T_f(\sigma)) = f^\ast \div(\sigma).
\]

\begin{theorem}\label{thm-difference-height-and-distance}
  Assume that $L$ and $f^\ast L$ are base-point-free. Suppose that $(L \cdot f(x)) >0$ for every
  $x\in \Ind(f)$.

  Then there exist strongly adelic metrics on $L$ and $f^\ast L$ and a family of constants
  $\{C_w\}_{w\in \Omega_K}$ such that, for any nonzero section $\sigma\in H^0(X,L)$, let
  $\sigma'=T_f(\sigma)$,

  \begin{enumerate}
    \item the map $(w\in \Omega_K^\infty) \mapsto C_w$ is $\nu$-dominated,
    \item the constants $C_w$ vanish for all but finitely many non-archimedean places
      $w\in \Omega_K^{\mathrm{NA}}$,
    \item for any $x\in (X\backslash\Ind(f))(K)$ and $f(x) \not\in \supp(\div(\sigma))$, we have
      \begin{align*}
        h_{(\overline{L},\sigma),w}(f(x)) - h_{(\overline{f^\ast L},\sigma'),w}(x) \leq \min \{0, \log \dist_w(x,\Ind(f)) +C_w\}.
      \end{align*}
  \end{enumerate}

\end{theorem}

\begin{proof}

  If $T_f(\sigma) (x) \not =0$ for some $x\in \Ind(f)$, then $\pi^\ast T_f(\sigma)$ does not vanish
  along $\pi^{-1}(x)$. Since $\pi^\ast f^\ast L - g^\ast L$ is effective, $g^\ast \sigma$ does not
  vanish along $\pi^{-1}(x)$. This implies $\sigma\cdot f(x) =0$ which contradicts $L\cdot f(x) >0$.
  Therefore, $\sigma'$ vanishes at $\Ind(f)$.

  Let $\Sigma = \{\tau_1,\ldots,\tau_s \}$ be a finite generating set of $H^0(X,L)$. For each
  $i=1,\ldots, s$, set $\tau'_i\coloneqq T_f(\tau_i)$. Then we extend $\{\tau'_1,\ldots,\tau'_s \}$
  to
  \begin{align*}
    \Sigma' = \{\tau_1',\ldots,\tau_s',\tau'_{s+1},\ldots,\tau'_t \}
  \end{align*}
  so that $\Sigma'$ generates $H^0(X,f^\ast L)$. Then $\Sigma$ and $\Sigma'$ define strongly adelic
  metric families on $L$ and $f^\ast L$. Let
  $$
    \varphi_w(x) = h_{(\overline{L},\sigma),w}(f(x)) - h_{(\overline{f^\ast L},\sigma'),w}(x).
  $$

  By definition, for $x\not\in \Ind(f)$ and $f(x)\not\in \supp(\div(\sigma))$, we have
  \begin{align*}
    \varphi_w(x) & = \log \max_{i=1,\ldots,s}\left| \frac{\tau_i(f(x))}{\sigma(f(x))}\right|_w - \log\max_{i=1,\ldots,t}\left| \frac{\tau'_i(x)}{\sigma'(x)}\right|_w\\
    & =\log\frac{\max_{i=1,\ldots,s}|\tau_i'(x)|_w}
    {\max_{i=1,\ldots,t}|\tau_i'(x)|_w}
  \end{align*}

  Note that $\tau_i'=\tau_i\circ f$ is globally defined and vanishes at $\Ind(f)$ for
  $i=1,\ldots,s$. It suffices to prove that, for each fixed $i_0$, the map
  $(w\in \Omega_K) \mapsto \log\sup_{x\in X(K)}\psi_{i_0,w}(x)$ is $\nu$-dominated, where
  \begin{align*}
    \psi_{i_0,w}(x) \coloneqq \frac{\min_{j=1,\ldots, t}\left|\tau'_{i_0}(x)/\tau'_j(x)\right|_w}{\dist_w(x,\Ind(f))}.
  \end{align*}

  Let $\cF$ be the sheaf defined by the divisors $\tau'_1,\ldots,\tau'_s$. Then its co-support is
  $\Ind(f)$. Use the fixed embedding $\iota: X\hookrightarrow \IP^N_K$. Choose a section $s_0$ of
  $\iota^\ast \cO_{\IP^{N}}(1)$ such that $\div(s_0)\cap \Ind(f)=\emptyset$. We consider
  $\cF \to \cF\otimes \iota^\ast \cO_{\IP^{N}}(a)$ for some $a\in \N$ large enough and identify each
  $\tau'_i$ with the restriction of a degree-$a$ homogeneous section
  $\rho_i\in H^0(\IP^N, \cO_{\IP^{N}}(a))$ on $\IP^N_K$.
  Then we may compute $\psi_{i_0,w}(x)$ on $\IP^N$ by $\rho_i$
  \begin{align*}
    \psi_{i_0,w}(x) \coloneqq \frac{\min_{j=1,\ldots, t} |\rho_{i_0}(x)/\rho_j(x) |_w}{\dist_w(\iota(x),\iota(\Ind(f)))}.
  \end{align*}

  We first consider the case where $w$ is non-archimedean. Fix an indeterminacy point
  $y=[y_0:\cdots:y_N]\in \iota(\Ind(f))$ and note that $\rho_{i_0}$ vanishes at $y$. Suppose that
  $y = [1:0:\cdots:0]$. Then we consider a neighborhood of $y$ with respect to $\dist_w$ which is
  $U_y= \{[1:x_1:\cdots:x_N], |x_i|_w <1\}$. For $x= [1:x_1:\cdots:x_N]\in U_y$, we write
  $\rho_{i_0}(x) = \sum_{1\leq|\alpha|\leq a} c_\alpha (x-y)^\alpha$, where $c_\alpha \in K$ and
  $(x-y)^\alpha = \prod_{i=1}^N (x_i-y_i)^{\alpha_i}$.
  Then for every $w\in \Omega_K^{\mathrm{NA}}$, one has
  \begin{align*}
    |\rho_{i_0}(x)|_w \leq \max_{1\leq|\alpha|\leq a} |c_\alpha|_w |(x-y)^\alpha|_w \leq C'_w(i_0,y)\dist_w(x,y),
  \end{align*}

  where
  $$
    C'_w(i_0,y) = \max_{1\leq|\alpha|\leq a} |c_\alpha|_w.
  $$
  Then we have
  $$
    \log\sup_{x\in X(K)}\psi_{i_0,w}(x) \leq \log\max_{y\in \Ind(f)}\{1,C'_w(i_0,y) \},
  $$
  and the right-hand side is zero for all but finitely many non-archimedean places $w$. Hence it is
  $\nu$-dominated.

  We now consider the case where $w$ is archimedean. Let $\pi\colon \sX \to \sB$ be a projective
  model, and let $\sL$ be a line bundle on $\sX$ such that $\sL_K = L$. The morphism $\iota$ extends
  to
  \[
    \widetilde \iota : \sX \to \IP^N_\sB.
  \]

  Recall that we view the set of archimedean places as a subset of $\sB(\IC)$, and that $\nu$ is the
  Monge--Amp\`ere measure associated with a metrized ample line bundle on $\sB$. The birational map
  $f\colon X\dasharrow X$ gives rise to $F\colon \sX\dasharrow \sX$. The sections $\rho_i$ induce
  $\sT_i$ for $i=1,\cdots, t$. Note that, for generic $w\in \sB(\IC)$,
  $F_w\colon \sX_w\dasharrow \sX_w$ is well defined and $\Ind(F) \supset {\Ind(F_w)}$ as well as
  $\Ind(F) \supset \overline{\Ind(f)}$.

  Let $\mathscr{F}$ be the ideal sheaf on $\IP^N_\sB$ defined by a section $\sT_{i_0}$ on
  $\IP^N_\sB$, for a fixed $i_0\in \{1,\ldots,s\}$. Let $\cG$ be the ideal sheaf on $\sX$ generated
  by $\sT_1,\ldots,\sT_t$. Set
  $V = \operatorname{Bs}(\sT_1,\ldots,\sT_s) = \overline{\operatorname{Bs}(\tau_1',\ldots,\tau_s')} \cup R$,
  where $R$ is the vertical component of $V$. Since $f^\ast L$ is globally generated by
  $\tau_1',\ldots,\tau_t'$, we note that $W = \operatorname{Bs}(\cG)$ is vertical and $W\subset R$.

  The \L ojasiewicz inequality applied to $\cF$ and $\cG$ yields constants $C''$ and $\alpha$ such
  that
  \begin{align*}
    g_{\cF}(p) \leq \log \dist_\sX(p, V) +C'',
  \end{align*}
  and
  \begin{align*}
    \alpha \log \dist_\sX(p, W) \leq g_{\cG}(p).
  \end{align*}
  Observe that, for $w\in \sB(\IC)$ and $p\in \sX_w(\IC)$,
  $$
    \log \dist_\sX(p, W) \geq \log\dist_{\sB}(\pi(p),\pi(W))+C''',
  $$
  for some constant $C'''$ independent of $w$, and
  \begin{align*}
    \log \dist_\sX(p, V) & \leq \log \dist_\sX(p, \overline{\Ind(f)} )\\
    & \leq \log \dist_{\sX_w}(p,\Ind(f_w)) =\log\dist_w(p,\Ind(f)).
  \end{align*}
  We conclude that
  \begin{align*}
    \log\psi_{i_0,w}(x) & = g_{\cF} - g_\cG - \log\dist_w(p,\Ind(f)) \\
    & \leq C''-\alpha C''' - \alpha \log\dist_{\sB}(w,\pi(W)).
  \end{align*}
  Then the map
  $(w\in \Omega_K^\infty) \mapsto C_w= C''-\alpha C''' - \alpha \log\dist_{\sB}(w,\pi(W))$ is
  $\nu$-dominated.
\end{proof}

\subsection{Lower bounds for periodic points}

\begin{theorem}\label{theorem-main-lower-bound}
  Let $f\colon X\dasharrow X$ be a birational self-map of a smooth projective surface defined over
  an algebraically closed field of characteristic $0$. Assume that the first dynamical degree
  $\lambda_{1}>1$. Then
  \begin{align*}
    \NIPer_n(f)\geq \lambda_{1}^n+o(\lambda_{1}^n).
  \end{align*}
\end{theorem}

We may assume that $f\colon X\dasharrow X$ is defined over $\IC$. If $f$ is birationally conjugate
to an automorphism, then the result follows from \cite{cantat2001autom},
\cite{bedford1993distribution}, and \cite{xie2014dynamique}. Thus it suffices to consider the case
when $f$ is strictly birational, i.e., there is no automorphism in its birational conjugacy class.
Then \cite{diller2001dynamics} shows that $X$ must be a rational surface. Again by
\cite{diller2001dynamics}, after finitely many blow-ups of $X$, a conjugate of $f$ satisfies
\begin{align*}
  \label{AS1} f^j(\Ind(f^{-1}))\cap \Ind(f) =\emptyset, \quad \text{for every $j\geq 0$}. \tag{AS1}
\end{align*}
Moreover, the first dynamical degree $\lambda_{1}$ is a simple eigenvalue for the action of $f^\ast$
on $H^{1,1}(X,\IR)$ and its eigenspace is generated by a nef class $\theta^+\in \Pic(X)_\IR$. Let
$\theta^{-1}\in \Pic(X)_\IR$ be the corresponding nef class of $f^{-1}$.

By \cite[Section~4]{bedford2005energy}, after finitely many blow-downs, one obtains a conjugate
satisfying \textup{(AS1)} as well as
\begin{align*}
  \label{AS2} (\theta^+\cdot f(x)) >0 ,\  \text{for $x\in \Ind(f)$}\quad\text{and}\quad (\theta^+\cdot f^{-1}(x)) >0,\ \text{for $x\in \Ind(f^{-1})$ } \tag{AS2}
\end{align*}

\begin{remark}
  Under these hypothesis, Jonsson--Reschke studied the dynamics of the base curves; see
  \cite[Section~4]{jonsson2018complex}. The set of base curve is
  $$
    \cC^+ =\{\text{irreducible curve $C$ on $X$}, (C\cdot\theta^+) =0\}.
  $$
  Set $\cC_{\text{per}} = \{C\in \cC^+, \text{$C$ is not $f^n$-exceptional for any $m\in \IN$}\}$.
  They showed that for any $C\in \cC_{\text{per}}$, we have $C\cap \Ind(f^n) =\empty$ for all
  $n\in \Z$ and $C$ does not intersect any $f^n$-exceptional curves for any $n\in \Z$.
\end{remark}

By collecting the coordinates and coefficients of $X$ and $f$, we may assume that such a birational
map $f\colon X\dasharrow X$ is defined over a field $K$ finitely generated over $\IQ$. Moreover, by
passing to a finite extension, we suppose the indeterminacy loci $\Ind(f)$ and $\Ind(f^{-1})$ are
also defined over $K$. We fix a model $(\sB,\overline\sH)$ of $K$ as in Section~\ref{section4-1},
and we have the space of places $\Omega_K = \Omega_K^{\mathrm{NA}}\sqcup\Omega_K^\infty$ and a
measure $\nu$ on $\Omega_K$.

\begin{theorem}\label{theorem-BD-condition-technical}
  Let $f\colon X\dasharrow X$ be the birational map considered above. For every point $y\in X(K)$
  whose forward orbit is well defined, we have
  \begin{align}\label{BD-condition}
    \sum_{j=0}^\infty \lambda_{1}^{-j}\log\dist_w(f^j(y),\Ind(f))> -\infty
  \end{align}
  for $\nu$-almost every place $w\in \Omega_K$.

  In particular, there exists an embedding $K\hookrightarrow \IC$ such that the base-changed map
  $f_\IC$ satisfies
  \begin{align}
    \sum_{j=0}^\infty \lambda_{1}^{-j}\log\dist(f_\IC^j(y),\Ind(f_\IC))> -\infty, \quad\text{for every $y\in \Ind(f_\IC^{-1})$}. \label{BD} \tag{BD}
  \end{align}
\end{theorem}

\begin{proof}
  If the Zariski closure of the orbit of $y$ is disjoint from $\Ind(f)$, then \eqref{BD-condition}
  holds for any embedding of $K$ into $\IC$. Therefore, by Lemma 4.9 of \cite{jonsson2018complex},
  we may assume that the Zariski closure of the orbit is an infinite set intersecting $\Ind(f)$, and
  that for $n\gg 0$, $f^n(y)$ does not belong to any of the curves in $\cC^+$, where $\cC^+$ is the
  set of irreducible curves $C$ on $X$ such that $(C\cdot \theta^+)=0$.

  By \cite{kawaguchi2008projective}, we may write $\theta^+ = A+[Z]$, where $A$ is ample and $Z$ is
  an effective divisor whose components all belong to $\cC^+$. Since $\theta^+$ is nef, for any
  strongly adelic metric $\varphi$ on $\theta^+$,
  Proposition~\ref{proposition-boundedness-effective-divisor-height} gives a constant $C_1$ such
  that
  \begin{align*}
    \frac{1}{\lambda_{1}^n} h_{\theta^+,\varphi}(f^n(y)) - h_{\theta^+,\varphi}(y) \geq - h_{\theta^+,\varphi}(y) +C_1
  \end{align*}
  for every $n\in\N$.

  We next use the decomposition of $\theta^+$ from Proposition 4.10 of \cite{jonsson2018complex}
  \begin{align*}
    \theta^+ = \sum_{i=1}^p r_iB_i + L + [D_{\text{per}}],
  \end{align*}
  such that
  \begin{enumerate}
    \item for each $i$, $r_i>0$; $B_i$ and $f^\ast B_i$ are base-point-free; and $(B_i\cdot f(x))>0$
      for every $x\in \Ind(f)$;
    \item $L\in \Pic_\IR(X)$ is nef;
    \item $D_{\text{per}}$ is an effective $\IR$-divisor supported on the curves in
      $\cC_{\text{per}}$.
  \end{enumerate}

  After passing to a finite extension of $K$, we may assume that all classes are defined over $K$.
  Let $B = \sum_{i=1}^p r_iB_i$. By Theorem~\ref{thm-difference-height-and-distance}, we obtain
  strongly adelic metrics on $B_i$ and $f^\ast B_i$, and sections $\sigma_{i,j}$ and
  $\sigma_{i,j}'=T_f(\sigma_{i,j})$, such that
  \begin{align*}
    h_{\overline{B_i},\sigma_{i,j} ,w}(f^{j+1}(y)) - h_{\overline{ f^\ast B_i},\sigma'_{i,j},w}(f^{j}(y)) \leq \log\dist_w(f^{j}(y),\Ind(f)) + C_w.
  \end{align*}

  Note that $C_w$ is $\nu$-dominated, which means that $\int_{w\in\Omega_K} |C_w|\,d\nu(w)<+\infty$.
  Therefore
  \begin{align*}
    & \sum_i r_i\int_{w\in \Omega_K} \sum_{j=0}^N \lambda_{1}^{-j}\log\dist_w(f^{j}(y),\Ind(f)) d\nu(w)\\
    \geq & \sum_i r_i\int_{w\in \Omega_K}\sum_{j=0}^N\lambda_{1}^{-j}\left( h_{\overline{B_i},\sigma_{i,j} ,w}(f^{j+1}(y)) - h_{\overline{f^\ast B_i},\sigma'_{i,j},w}(f^{j}(y))\right)d\nu(w) -\int_{w\in\Omega_K} |C_w|\,d\nu(w)\\
    = & \sum_i r_i\sum_{j=0}^N\lambda_{1}^{-j}\left( h_{\overline{B_i}}(f^{j+1}(y)) - h_{\overline{f^\ast B_i}}(f^{j}(y))\right)-\int_{w\in\Omega_K} |C_w|\,d\nu(w).
  \end{align*}

  Equipping $L$, $f^\ast L$, $D_{per}$, and $f^\ast D_{per}$ with adelic metrics, we obtain adelic
  metrics $\varphi_1$ on $\theta^+$ and $\varphi_2$ on $f^\ast \theta^+ = \lambda_{1} \theta^+$.

  First, we have
  \begin{align*}
    & \sum_{j=0}^N\lambda_{1}^{-j} \left( h_{\theta^+,\varphi_1}(f^{j+1}(y))- h_{f^\ast \theta^+,\varphi_2}(f^{j}(y))\right) \\
    = & \sum_{j=0}^N\lambda_{1}^{-j} \left( h_{ \theta^+,\varphi_1}(f^{j+1}(y)) - h_{\lambda_{1} \theta^+,\lambda_{1} \varphi_1}(f^{j+1}(y)) \right)\\
    & \qquad\qquad\qquad\qquad\qquad + \sum_{j=0}^N\lambda_{1}^{-j} \left( h_{f^\ast \theta^+,\lambda_{1}\varphi_1}(f^{j+1}(y)) - h_{f^\ast \theta^+,\varphi_2}(f^{j}(y)) \right)\\
    & = -\dist(\lambda_1 \varphi_1,\varphi_2) + \frac{1}{\lambda_{1}^{N+1}} h_{\theta^+,\varphi_1}(f^{N+1}(y)) - h_{\theta^+,\varphi_1}(y)\\
    & \geq -\dist(\lambda_1 \varphi_1,\varphi_2) +C_1- h_{\theta^+,\varphi_1}(y) >-\infty.
  \end{align*}

  Second, by Proposition~\ref{prop-difference-height-nef-divisor}, we have constant $C_2, C_3$
  independent on $j$ such that
  \begin{align*}
    & \left( h_{\overline{ L}}(f^{j+1}(y)) - h_{\overline{ f^\ast L}}(f^{j}(y))\right) \geq C_2.
  \end{align*}
  \begin{lemma}\label{Lemma-difference-height-periodic-curve}
    For any strongly adelic metrics on $\cO(C)$ and $\cO( f^\ast C)$, there is a constant $M$ such
    that for $x\not\in \Ind(f)$ with $f(x)\not\in C$,
    \begin{align*}
      \left|h_{\overline{C}}(f(x)) - h_{\overline{f^\ast C}}(x)\right| \leq M.
    \end{align*}

  \end{lemma}
  So there exists a constant $C_3$
  \[
    \left| h_{\overline{ [D_{\text{per}}]}}(f^{j+1}(y)) - h_{\overline{f^\ast [D_{\text{per}}]}}(f^{j}(y))\right| \leq C_3 .
  \]
  In conclusion, we obtain the following uniform estimate:
  \begin{align*}
    \int_{w\in \Omega_K} & \sum_{j=0}^N \lambda_{1}^{-j}\log\dist_w(f^{j}(y),\Ind(f)) d\nu(w) \geq C_4.
  \end{align*}
  where
  \[
    C_4 = -\frac{1}{\sum r_i}\int_{w\in\Omega_K} |C_w|\,d\nu(w)-\dist(\lambda_1 \varphi_1,\varphi_2) +C_1- h_{\theta^+,\varphi_1}(y) -C_2 -C_3.
  \]
  This in turn implies
  \begin{align*}
    \int_{w\in \Omega_K} & \sum_{j=0}^{+\infty} \lambda_{1}^{-j}\log\dist_w(f^{j}(y),\Ind(f)) d\nu(w) > -\infty.
  \end{align*}
  One concludes that for $\nu$-almost every $w\in \Omega_K^\infty$ and for every non-archimedean
  place $w\in \Omega_K^{\mathrm{NA}}$, we have
  \begin{align*}
    \sum_{j=0}^{+\infty} \lambda_{1}^{-j}\log\dist_w(f^{j}(y),\Ind(f)) > -\infty.
  \end{align*}

  Since the set $\Ind(f^{-1})$ is finite, condition~\ref{BD} holds for $\nu$-almost every
  archimedean place. Hence such an embedding exists.
\end{proof}

\begin{proof}[Proof of Lemma~\ref{Lemma-difference-height-periodic-curve}]
  We take the minimal resolution $g\circ\pi^{-1}$ of $f$. By \cite{jonsson2018complex}, the pullback
  $C' \coloneqq f^\ast(C)$ belongs to $\cC_{\text{per}}$ and
  $C'\cap \Ind(f) = C\cap \Ind(f^{-1}) = \emptyset$. Hence $g^\ast C = \pi^\ast C'$. For
  $x\not\in \Ind(f)$, there is a unique preimage $y = \pi^{-1}(x)$. Then
  \begin{align*}
    h_{\overline{[C]}}(f(x)) = h_{\overline{[C]}}(g(y)) = h_{g^\ast \overline{[C]}}(y)= h_{\pi^\ast\overline{[C']}}(y) +O(1) = h_{\overline{[C']}}(x)+O(1),
  \end{align*}
  which completes the proof.
\end{proof}

\begin{proof}[End of proof of Theorem~\ref{theorem-main-lower-bound}]

  Assuming first $f$ satisfies (AS1) and (AS2). By Theorem~\ref{theorem-BD-condition-technical},
  there exists an archimedean place $w\in \Omega_K^\infty$ such that for every $y\in \Ind(f^{-1})$,
  the condition \ref{BD-condition} holds. Thus there exists an embedding $K\hookrightarrow\IC$ such
  that, after the base change, $f_\IC$ satisfies the BD condition \cite{bedford2005energy}. By the
  main theorem of \cite{bedford2005energy} and Theorem~2 in \cite{dujardin2006laminar}, there is an
  $f_\IC$-invariant measure $\mu$ on $X$ that puts no mass on any proper algebraic subset. Further,
  we have that the saddle periodic points of $f_\IC$ are equidistributed with respect to $\mu$. Let
  $U$ be any Zariski open subset of $X$. Since $\mu$ has full mass on $U$, we have
  \begin{align*}
    \lim_{n\to \infty}\frac{1}{\lambda_{1}^n}\#\{x\in U, \text{$x$ is a saddle periodic point under $f_\IC^n$}\} = 1.
  \end{align*}

  Therefore, for any conjugate $f'$ of $f$, we consider the same embedding as above, then
  \[
    \lim_{n\to \infty}\frac{1}{\lambda_{1}^n}\#\{x\in U, \text{$x$ is a saddle periodic point under $f'_\IC$} \}= 1
  \]
  Since all saddle periodic points are isolated, we concludes that
  \begin{align*}
    \NIPer_n(f') \geq \lambda_{1}^n+o(\lambda_{1}^n).
  \end{align*}

\end{proof}
Combining Theorem~\ref{upper-bound-surface} and Theorem~\ref{theorem-main-lower-bound}, we
immediately obtain Theorem~\ref{theoremB}.

\appendix

\section{A proof of Theorem~\ref{theoremA} on normal projective variety}

The ideas leading to the next theorem were suggested by Chat-GPT. The author subsequently wrote the
details of the proof.
\begin{theorem}
  Let $X$ be a normal projective variety of dimension $m$ over an algebraically
  closed field $K$, and let $L$ be an ample line bundle on $X$. There exists
  a constant $C_{X,L}>0$ such that, for every dominant rational self-map $f$
  of $X$ and every $n\geq1$,
  \[
    \NIPerM_n(f)\leq C_{X,L}\max_{0\leq p\leq m}\deg_{L,p}(f^n).
  \]
\end{theorem}

\begin{proof}
  Fix $n\geq1$ and write $\Gamma=\Gamma(f^n)$ and
  $Z=\Delta\times_{X\times X}\Gamma$. For $m=0$ both the count and the
  degree equal one. If $m>0$ and $f^n=\mathrm{id}$, then $Z=\Delta$ has no
  isolated points. We may thus assume $m\geq1$ and $\Gamma\not\subseteq\Delta$.

  \smallskip\noindent\textbf{Step 1.}
  Let $\pi: \widetilde{X\times X} =\operatorname{Bl}_\Delta(X\times X)\to X\times X$ have
  exceptional divisor $E$, and let $\widetilde\Gamma$ be the strict transform of $\Gamma$. Let
  \(p_1, p_2\) the projections \(X\times X\to X\).
  Choose $M\geq1$, depending only on $(X,L)$, so that
  \[
    L_1=M (p_1^*L+p_2^*L) \quad\text{and}\quad L_2=\pi^*L_1-E
  \]
  are very ample on $X\times X$ and $\widetilde{X\times X}$, respectively.

  We consider the scheme
  $$
    R=E\cap_{ \widetilde{X\times X}} \widetilde\Gamma.
  $$
  Note that the intersection \(\widetilde\Gamma\) with \(E\) is proper, so
  $[\widetilde\Gamma]\cdot c_1(E)=[R]$ is effective. Using \(\pi^\ast L_1 - L_2 = E\), we have
  \[
    (\pi^\ast L_1)^m\cdot[\widetilde\Gamma]-{L_2}^m\cdot[\widetilde\Gamma]
    =\left( \sum_{j=0}^{m-1} (\pi^\ast L_1)^{m-1-j} L_2^j \right)\cdot( \pi^\ast L_1 - L_2)\cdot[\widetilde\Gamma].
  \]
  Both $\pi^\ast L_1$ and $L_2$ are nef and. So the right hand side is non-negative. The projection
  formula gives
  \[
    L_1^m \cdot [\Gamma]= (\pi^\ast L_1)^m\cdot[\widetilde\Gamma]\geq {L_2}^m\cdot[\widetilde\Gamma].
  \]
  Thus we conclude that

  \[
    \operatorname{Vol}_{L_2}(\widetilde\Gamma) \leq \operatorname{Vol}_{L_1}(\Gamma).
  \]

  Since \(L_2 -E\) is still nef, we have
  \[
    \Vol_{L_2}(R) = L_2^{m-1}\cdot[R] = L_2^{m-1}\cdot E \cdot [\widetilde\Gamma] \leq L_2^{m-1}\cdot [\widetilde\Gamma] = \Vol_{L_1}(\Gamma).
  \]

  \smallskip\noindent\textbf{Step 2.}
  The strict transform \(\widetilde \Gamma\) is canonical isomorphic to
  $\operatorname{Bl}_{Z}\Gamma$ the blow-up of \(\Gamma\) along a subscheme \(Z\) by the blow-up
  closure lemma \cite[\href{https://stacks.math.columbia.edu/tag/080E}{Tag 080E}]{stacks-project}.
  And $R$ is its exceptional divisor.

  For each isolated point $P=(x,x)$ of $Z$, put $Z_x=\operatorname{Spec}(S/I)$. As $Z_x$ is open and
  closed in $Z$, so is
  $R_x=R\times_Z Z_x$ in $R$. In particular,
  $[R]-\sum_{P\in Z^{\mathrm{iso}}}[R_x]$ is effective.

  Observe that \(R_x\) is the exceptional divisor of the blowing-up of \(\Gamma\) at \(P = (x,x)\).
  Over $Z_x$ the line bundle $L_1$ is trivial, so we have \(L_2|_{R_x}\simeq\mathcal O_{R_x}(1)\).
  Put \(h = c_1(L_2)\).

  The multiplicity of \(\Gamma\) along \(P\), in
  Section~\ref{section-Multiplicity along subvariety}, is
  \[
    e_P(\Gamma) = c_1(\cO_{R_x}(1))\cdot [R_x] = h^{m-1}\cdot[R_x].
  \]

  Intersecting the remaining effective cycle with $h^{m-1}$ yields
  \[
    \sum_{P=(x,x)\in Z^{\mathrm{iso}}} e_x(f^n)
    \leq h^{m-1}\cdot[R]\leq \operatorname{Vol}_{L_1}(\Gamma).
  \]

  \smallskip\noindent\textbf{Step 3.} We claim that there exists a constant
  $c_X$ depending only on $X$ such that, for every $x\in\IPer_n(f)(K)$,
  \[
    l_x(f^n)\leq c_X e_x(f^n).
  \]
  Indeed, for such a point $x$ and put $P=(x,x)$. Since the length-$n$ orbit of $x$
  is defined, the first projection identifies $\mathcal O_{\Gamma,P}$ with
  $S=\mathcal O_{X,x}$. Let $I\subset S$ be the ideal of $Z$ under this
  identification. The point $P$ is isolated in $Z$, so $I$ is
  $\mathfrak m_x$-primary. By definition,
  \[
    l_x(f^n)=\ell_S(S/I),\qquad e_x(f^n)=e(I,S).
  \]

  On an affine chart $U=\operatorname{Spec}A$ of $X$, apply
  \cite[Remark~2.6]{klein2019lech} to $A$ as a module over itself. The finite type $K$-domain $A$ is locally equidimensional and is a quotient of a polynomial ring over $K$. Its residue fields are infinite, and it satisfies the uniform Artin--Rees property \cite[Appendix~A]{klein2019lech}. Thus there is a constant $c_U$ such that $\ell_{A_{\mathfrak p}}(A_{\mathfrak p}/J)\leq c_U e(J,A_{\mathfrak p})$
  for every $\mathfrak p\in\operatorname{Spec}A$ and every
  $\mathfrak p A_{\mathfrak p}$-primary ideal $J$. Taking the maximum of
  the $c_U$ over a finite affine cover and specializing to $S$ and $I$
  proves the claim.

  \smallskip\noindent\textbf{Conclusion.}
  Step~2 and Step~3 gives that
  \begin{align*}
    \NIPerM_n(f) = \sum_{x\in \IPer_n(f)}l_x(f^n)
    & \leq c_X\sum_{x\in \IPer_n(f)} e_x(f^n) \leq c_X\sum_{P=(x,x)\in Z^{\mathrm{iso}}}e_x(f^n)
    \\
    & \leq c_X\,\operatorname{Vol}_{L_1}(\Gamma) \leq c_X(2M)^m\max_{0\leq p\leq m}\deg_{L,p}(f^n).
  \end{align*}
  The last inequality is from Section~\ref{section-upper-bound-rational-map}. This proves the
  assertion with $C_{X,L}=c_X(2M)^m$.
\end{proof}

\bibliographystyle{alpha}
\bibliography{ag}

\end{document}